\documentclass[a4paper,oneside,11pt,reqno]{amsart}
\usepackage[backend=bibtex,style=ieee,citestyle=numeric-comp,giveninits=true,url=false,eprint=false,doi=false]{biblatex}
\usepackage{hyperref}
\usepackage[foot]{amsaddr}
\usepackage{amsmath}
\usepackage{amsthm}
\usepackage[english]{babel}
\usepackage{xurl}
\usepackage{graphicx}

\usepackage{physics}

\usepackage{subcaption}
\usepackage{todonotes}
\usepackage{algorithm}
\usepackage{algpseudocode}

\usepackage{empheq}

\usepackage{xcolor}

\usepackage[T1]{fontenc}
\usepackage[utf8]{inputenc}
\usepackage{textcomp} 
\usepackage{amssymb}
\usepackage{mathtools}

\usepackage[scr=boondox]{mathalpha}
\newcommand{\symbb}[1]{\mathbb{#1}}
\newcommand{\symcal}[1]{\mathcal{#1}}
\newcommand{\symscr}[1]{\mathscr{#1}}
\newcommand{\symfrak}[1]{\mathfrak{#1}}
\renewcommand{\smallsetminus}{\setminus}

\calclayout

\newcommand{\domain}{D}
\DeclareMathOperator{\tDom}{[0,T)}
\newcommand{\tInt}{\int_{0}^{T}}
\newcommand{\soln}{u}
\newcommand{\testSoln}{\varphi}
\newcommand{\testEnt}{\xi}
\newcommand{\entropy}{\eta}
\newcommand{\entflux}{q}
\DeclareMathOperator{\residual}{\symcal{R}}
\DeclareMathOperator{\resfun}{\symscr{r}}
\DeclareMathOperator{\Loss}{\symcal{L}}
\DeclareMathOperator{\dLoss}{\symfrak{L}}
\newcommand{\activation}{\mathop\sigma} 
\DeclareMathOperator{\sgn}{sign}
\DeclareMathOperator{\relerr}{\symcal{E}_r}
\newcommand{\collocPoints}{S}

\theoremstyle{definition}
\newtheorem{definition}{Definition}
\theoremstyle{Theorem}
\newtheorem{remark}{Remark}

\begin{document}
\title[Weak PINNs for hyperbolic conservation laws]{Improving Weak PINNs for Hyperbolic Conservation Laws:\\ Dual Norm Computation, Boundary Conditions and Systems}
\author{Aidan Chaumet}
\email[Aidan Chaumet]{chaumet@mathematik.tu-darmstadt.de}
\author{Jan Giesselmann}
\email[Jan Giesselmann]{giesselmann@mathematik.tu-darmstadt.de}

\begin{abstract}
	We consider the approximation of entropy solutions of nonlinear hyperbolic conservation laws using neural networks. We provide explicit computations that highlight why classical PINNs will not work for discontinuous solutions to nonlinear hyperbolic conservation laws and show that weak (dual) norms of the PDE residual should be used in the loss functional. This approach has been termed “weak PINNs” recently. We suggest some modifications to weak PINNs that make their training easier, which leads to smaller errors with less training, as shown by numerical experiments. Additionally, we extend wPINNs to scalar conservation laws with weak boundary data and to systems of hyperbolic conservation laws. We perform numerical experiments in order to assess the accuracy and efficiency of the extended method.
\end{abstract}
\date{\today}

\maketitle
	\date{\today}
	
\section{Introduction}
Machine learning is an effective approach for a variety of challenging applications such as image generation, computer vision or natural language processing. Their ubiquitous success has triggered mathematical studies into the use of neural networks for approximating solutions to partial differential equations. Physics informed neural networks (PINNs) as introduced in \cite{raissi2019physics} are a popular approach for approximating PDEs. Classical PINNs leverage automatic differentiation for neural networks to compute the PDE residual and minimize its $L^2$ norm in order to approximate solutions. For problems with smooth solutions, error estimates such as \cite{deRyck_2022,mishra_L2fail} explain the success of PINNs. In contrary, empirical evidence suggests that classical PINNs do not work for discontinuous solutions \cite{Lye_2020} that naturally occur as solutions to nonlinear hyperbolic conservation laws. 

Hyperbolic conservation laws arise in many models in continuum physics. Examples include the Euler- and shallow water equations. The development of numerical schemes for hyperbolic conservation laws in up to three space dimensions is quite mature and PINNs are non-competitive for many classical applications. However, there are new challenges that require schemes which are efficient in high space dimensions, where mesh-based methods are unfeasible, e.g. the approximation of correlation measures that appear in the definition of statistical solutions \cite{Fjordholm2017}.

The first goal of this manuscript is to explain the failure of PINNs for discontinuous solutions of nonlinear hyperbolic conservation laws using explicit computations. This establishes that the empirical observations are not due to practical difficulties, such as insufficiently large networks or lack of training. Instead, we show that classical PINNs are fundamentally unable to approximate discontinuous solutions of nonlinear hyperbolic conservation laws due to their loss functional. Several elementary, but to the best of our knowledge new, computations show that good continuous approximations to discontinuous solutions do not produce residuals which are small in the $L^p$ norm for any $1 \leq p \leq \infty$. Thus, minimizing the $L^p$ norm of the PDE residual over the set of functions represented by a neural network cannot give good approximations.

Our computations show that instead, for good approximations, weak norms of the PDE residual are small. This is analogous to replacing the strong formulation of PDEs with their weak formulation to describe discontinuous solutions. Examples of PINNs using weak norms of PDE residuals to successfully approximate discontinuous solutions are variational PINNs \cite{variPinn} or weak PINNs \cite{wPINNS} (wPINNs). Both of these approaches approximate weak norms of the PDE residual by invoking the $\sup$-based definition of weak norms and maximizing over test functions. This gives a saddle-point problem during training. The loss is minimized with respect to the approximate solution and maximized with respect to the test function giving the weak norm of the PDE residual.
Training this type of PINNs is difficult, especially because the maximization problem may get stuck in bad local maxima. This also limits their attractiveness for concrete applications.

This motivates the second goal of this paper. In order to make training wPINNs easier, we approach the weak norm estimation differently. Instead of using the $\sup$-based definition of weak norms, we identify the maximizer of the dual pairing as the solution to a family of dual elliptic problems. Then, the solutions to the elliptic problems are identified as maximizers of a concave functional. We use a neural network to represent the solution and train it with gradient ascent. This is similar to the Deep Ritz Method for solving the Poisson problem \cite{E2018} and accelerates the learning of the dual problem, leading to overall more stable and effective learning of  wPINNs.

For many practical applications, further generalizations of wPINNs are required. We discuss extensions to weak boundary conditions and systems of hyperbolic conservation laws. The original wPINN framework prescribes Dirichlet boundary data for the solution on the entire boundary. However, prescribing Dirichlet boundary data is only possible on the (solution-dependent) inflow part of the boundary, so prescribing the correct boundary data requires prior knowledge of the solution. Secondly, the loss functional of original wPINNs is based on the Kruzhkov entropies, which is convenient because this formulation ensures one gets the unique weak entropy solution. However, the Kruzhkov entropies are only available for scalar hyperbolic conservation laws.

In light of this, our third goal is to extend wPINNs to incorporate weak boundary data understood in the sense of \cite{kondo2001measure}. This approach is valid for scalar conservation laws in multiple space dimensions. We replace the Dirichlet-based boundary contribution of the loss by a suitable entropy-based inequality on the boundary of the domain.

During preliminary numerical experiments using wPINNs with weak boundary conditions we observe that long time domains may lead to poor approximate solutions. Sometimes, wPINNs choose to make a large error at early times, in order to remove discontinuities from the domain and then minimize residuals for a smooth function subsequently. While leading to overall worse solutions, this may seem advantageous because the loss is only minimized locally in the neural network's parameter space. This is inherently tied to the fact that PINNs do not respect causality. We present a time-slicing approach that enforces causal relationships between the time slices, but not within a single time slice. This resolves the encountered difficulties.

Our last goal is to extend wPINNs to systems of hyperbolic conservation laws. This means replacing the  Kruzhkov entropy residual in the original wPINN loss by the PDE residual measured using a suitable weak norm. We do this in the scalar case already, as it seems to speed up training. Additionally, one has to ensure that one obtains the unique entropy solution. We fix a single strictly convex entropy-entropy flux pair and define an entropy residual that measures the violation of the corresponding entropy inequality. This is justified, because for scalar conservation laws with convex flux, satisfying a single entropy condition is equivalent to satisfying all entropy conditions. Further, for systems, there often is only a single physically motivated entropy. Our modified loss then naturally extends to systems of hyperbolic conservation laws by summing over the weak norms of the PDE residual of each individual equation.

Our manuscript is structured as follows. In section \ref{sec:lit_review} we give some background on various related machine learning methods used for solving PDEs. In section \ref{sec:hypCons} we provide new insights into why classical PINNs fail for discontinuous solutions of nonlinear hyperbolic conservation laws and why wPINNs are a suitable framework for their treatment. Then, section \ref{sec:weakNorms} presents a novel approach for approximating weak norms using neural networks, in order to make wPINNs more efficient. Further, we extend wPINNs to weak boundary conditions for scalar conservation laws and extend them to systems of hyperbolic conservation laws.
In section \ref{sec:numResults} we describe our training algorithm and then show numerical experiments comparing the performance of wPINNs as described in \cite{wPINNS} against our modified version, using the same examples as in \cite{wPINNS}. For some basic examples, original wPINNs and our modified wPINNs give nearly identical results. However, in the more challenging tests, the modified wPINNs clearly outperform the original wPINNs. Finally, in section \ref{sec:extensions}, we demonstrate the effectiveness of our extensions to weak boundary conditions and systems on some examples. In Appendix \ref{apx:time_dependent} we provide an additional, more general, computation explaining the failure of PINNs. In Appendix \ref{apx:implementation} we give further details on the implementation of our modified wPINNs in order to make results more transparent and reproducible.
	\section{Background on PINNs}
\label{sec:lit_review}

Originally, PINNs go back to \cite{PinnOrig1,PinnOrig2,PinnOrig3}, but were popularized more recently in \cite{raissi2019physics}. They are an unsupervised learning method, not requiring labeled data from observations or numerical schemes. However, when data is available it can be readily incorporated. Interest in PINNs has grown rapidly due to their versatility. One of their main benefits is their meshless nature, which may make them suitable for high-dimensional problems. Some recent examples of these methods include \cite{SirignanoDGM}, introducing the Deep Ritz Method, a PINN framework for elliptic problems, \cite{Ming_2021}, presenting advancements in enforcing essential boundary conditions for the aforementioned framework, \cite{fractionalADE}, wherein PINNs are applied to fractional advection-diffusion equations and \cite{Hu_2022}, proposing a domain decomposition framework for PINNs for multiscale problems and to improve parallelization during training. 

In their most frequently used form, PINNS are based on minimizing the $L^2$ norm of the residual associated with the PDE, admitting interpretation of PINNs as a weighted least-squares point-collocation method \cite{bochev2009least,PATEL2022110754}. There has been tremendous progress in the convergence analysis of PINNs \cite{Minakowski_2023,shin_2020,deRyck_2022}, that is based on three key principles \cite{wPINNS}:
\begin{enumerate}
	\item \textit{Regularity} of the solutions to the underlying PDEs [\ldots] which can be leveraged into proving that PDE residuals can be made arbitrarily small.
	\item \textit{Coercivity} (or stability) of the underlying PDEs, which ensures that the total error may be estimated in terms of residuals. For nonlinear PDEs, the constants in these coercivity estimates often depend on the regularity of the underlying solutions.
	\item \textit{Quadrature error} bounds for estimating the so-called generalization gap between the continuous and discrete versions of the PDE residual \cite{DERYCK2021732}.
\end{enumerate}
As mentioned in \cite{wPINNS}, this analysis is not applicable to weak solutions, because of their reduced regularity. In \cite{mishra_L2fail}, it is indeed observed in numerical experiments that classical PINNs fail for discontinuous solutions to nonlinear hyperbolic conservation laws. We explain this in section \ref{sec:hypCons}.
Because PINNs can be understood as a point-collocation method, they share some basic ideas with finite-element based point-collocation methods. For these, it has been observed that it is crucial which norm is used for minimization. In \cite{Guermond_Lp}, it is shown that for advection-reaction problems {with ill-posed boundary conditions}, $L^2$ minimization of PDE residuals does not recover the unique viscosity solution {in the sense of \cite{BLN_boundary}} to the problem. Instead, \cite{Guermond2008,Guermond_L1} show that for linear first-order PDEs in one dimension with ill-posed boundary conditions, minimizing residuals with respect to the $L^1$ norm in the space of piecewise linear functions for given boundary data is a fast and effective strategy that recovers the viscosity solution. In \cite{Guermond_L1}, the $L^1$ minimization technique is then extended to (some) non-linear PDEs in one dimension, such as the Hamilton-Jacobi equations. This is shown to be successful when the exact solution is continuous. However, \cite{Guermond_L1} provides an example where $L^1$ minimization fails to recover the correct discontinuous viscosity solution for the inviscid Burgers equation.

An established strategy to approach this problem is viscous regularization \cite{Lavery_burgers}.  	
Viscous regularization has been applied to PINNs too \cite{raissi2019physics}. However, this raises the question how much viscosity one needs to use to obtain solutions that can be learned effectively, yet are still close to the original solution. 

The choice of norm is important for PINNs as well and there has been recent research into suitable choices. In \cite{L2_HJB}, it is suggested that for high-dimensional Hamilton-Jacobi-Bellmann equations, one should choose an $L^p$-based loss with sufficiently large $p$, and a procedure for training with respect to an \mbox{$L^\infty$-based} loss is outlined. Approximating solutions with reduced regularity has led to multiple variants of the PINN methodology, such as control volume PINNs \cite{PATEL2022110754}, variational PINNs \cite{variPinn}, the Deep Ritz Method \cite{E2018} and weak PINNs \cite{wPINNS}.	
For control volume PINNs, the conservation law is expressed in integral form over test volumes, giving a formulation reminiscent of finite volume methods. The underlying idea of the other methods can be understood as minimizing the residual not in the $L^2$ norm, but in some weak (dual) norm.  While the Deep Ritz Method is specifically tailored to elliptic problems by using the fact that their solutions are minimizers of energy functionals, variational PINNs and weak PINNs address more general PDEs. 
Both of these methods measure the residual in a dual norm, but the key difference between the two is the way in which the test functions are represented. While variational PINNs use a mesh-based approach to represent test functions, weak PINNs represent the test functions by neural networks, thereby retaining the meshless nature of the method. The concept of weak PINNs was introduced on the prototypical example of scalar hyperbolic conservation laws in \cite{wPINNS}. The authors provide a detailed convergence analysis as well as numerical experiments demonstrating the capability to approximate discontinuous solutions. Another approach to solving problems with discontinuous solutions is outlined in \cite{liu2022discontinuity}, where the residual is minimized in a weighted $L^2$ norm with an adaptively chosen weight, chosen such that residuals close to steep gradients are weighted less.
\section{Explaining the failure of classical PINNs}
\label{sec:hypCons}
In this section, we present our analytical arguments that explain why classical PINNs cannot approximate discontinuous solutions to nonlinear hyperbolic conservation laws well. We begin by recalling some basic definitions involving scalar conservation laws and classical PINNs before showing that continuous approximations of discontinuous solutions have large $L^p$ norms of PDE residuals and thus will not be learned through minimizing the neural network loss. Then we motivate that weak norms do not have this issue.
\subsection{Scalar Conservation Laws}	
We consider the scalar conservation law given by
\begin{empheq}[left=\empheqlbrace]{align}
\soln_t + f(\soln)_x &= 0 &&  \text{in } \tDom \times \symbb{R}, \label{eqn:cLaw}\\
\soln(0,\cdot) &= \soln_0&&  \text{on } \symbb{R} \label{eqn:iCond},
\end{empheq}
where $\soln$ is the \textit{conserved quantity}, $f \in C^2(\symbb{R})$ is the \textit{flux function} and $u_0 \in L^\infty(\symbb{R})$ is the initial data. 
In the following, we always assume that $f \in C^2(\symbb{R})$ and $f$ is strictly convex. 

Oftentimes one wishes to work not on all of $\tDom \times \symbb{R}$, but rather restrict oneself to a compact interval $\domain \subset \symbb{R}$ in space. In this case one must supply boundary conditions, which need to be understood in a suitable weak sense, see \cite{BLN_boundary,kondo2001measure} and references therein. However, for initial conditions that are constant everywhere except on some compact interval in space, and boundary conditions that are constant in time, one knows that waves propagate with some maximum speed, such that the solution is constant outside some trapezoidal region in space-time.	 
For now, we limit our discussion to Dirichlet boundary conditions on space intervals $\domain$ such that the left- and right boundary are outside of this trapezoid for the given time interval $\tDom$. We denote the boundary conditions by $u(t,x) = g(t,x)$ for  $(t,x)\in (0,T) \times \partial \domain$, with $g$ satisfying the compatibility condition $u_0(x) = g(0,x)$ for $x \in \partial\domain$. Later, we extend wPINNs to weak boundary conditions.

Even for very simple flux functions and smooth initial conditions, one may obtain discontinuous solutions after finite time. Hence, one has to work with weak solutions:
\begin{definition}[Weak Solution]
	We call $\soln \in L^\infty([0,T) \times\symbb{R})$ a weak solution of problem \eqref{eqn:cLaw}-\eqref{eqn:iCond} for some initial datum $\soln_0 \in L^\infty(\symbb{R})$, if it satisfies  
	\begin{equation}
	\int_{[0,T)} \int_{\symbb{R}} \soln \testSoln_t + f(\soln) \testSoln_x \dd{x} \dd{t}
	+ \int_{\symbb{R}}\soln_0 \testSoln(0,\cdot) \dd{x} = 0,
	\label{eqn:wcLaw}
	\end{equation}
	for any test function $\testSoln \in C_c^1([0,T) \times\symbb{R})$.
	\label{def:weakSoln}
\end{definition}
Weak solutions are, in general, not unique. One strategy to establish uniqueness is by imposing \textit{entropy admissibility conditions} on weak solutions \cite{dafermos}.
We define entropies $\entropy$ and corresponding entropy fluxes $\entflux$ as follows.

\begin{definition}[Entropy -- Entropy Flux Pair]
	Let $\entropy,\entflux \in C^1(\symbb{R})$.  We say that $(\entropy,\entflux)$ are a \textit{convex entropy-entropy flux pair}, if $\entropy$ is strictly convex and $\entflux' = f'\entropy'$ almost everywhere.
\end{definition} 
Using the notion of entropy-entropy flux pairs one may define \textit{entropy-admissible} solutions to problem \eqref{eqn:cLaw}-\eqref{eqn:iCond}.
\begin{definition}[Entropy Admissibility]
	We call $\soln$ an \textit{entropy-admissible} solution to the Cauchy problem  \eqref{eqn:cLaw}-\eqref{eqn:iCond} if it is a weak solution in the sense of Definition \ref{def:weakSoln} and for any convex entropy-entropy flux pair, it satisfies 
	\begin{equation}
	\entropy(u)_t + \entflux(u)_x \leq 0 \quad \text{in } \tDom \times \symbb{R}
	\label{eqn:entDiss}
	\end{equation}
	in a distributional sense.
\end{definition}
For a strictly convex flux function $f$, 
weak solutions $\soln$ satisfying equation \eqref{eqn:entDiss} for a single strictly convex entropy-entropy flux pair are entropy-admissible \cite{panov1994uniqueness,de2004minimal}. Further, entropy-admissible solutions are unique.
\subsection{Classical PINNs}

In the following, we summarize the classical PINNs approach. We discuss fully connected feed-forward (FCFF) neural networks with hidden layers of uniform width to approximate solutions to scalar conservation laws. Note that there are a multitude of more sophisticated network architectures, as studied in other works, e.g. \cite{SirignanoDGM,E2018}. Mathematically we define such neural networks as follows:
\begin{definition}[Fully Connected NN]
	\label{def:fcNN}
	Let $l,w \in \symbb{N}$. Let $d\in\symbb{N}$ be the input dimension. Let $\activation:\symbb{R} \rightarrow \symbb{R}$ be a smooth, non-linear function. We call $\activation$ the \emph{activation function}. 
	
	Set $L^{[0]}(x) \coloneqq W^{[0]}x + b^{[0]}$ with $W^{[0]} \in \symbb{R}^{w\times d}$ and $b^{[0]} \in \symbb{R}^w$. Further, set $L^{[i]}(x) \coloneqq W^{[i]}x + b^{[i]}$ with $W^{[i]} \in \symbb{R}^{w\times w}$ and $b^{[i]} \in \symbb{R}^w$ for $i = 1,\dotsc,l-1$. Finally, set $L^{[l]}(x) \coloneqq W^{[l]}x + b^{[l]}$ with $W^{[l]} \in \symbb{R}^{1\times w}$ and $b^{[l]} \in \symbb{R}$. We call the matrices $W^{[i]}, i = 0,\dotsc,l$ the \emph{weights} and $b^{[i]}$ the \emph{biases} of a neural network. The affine linear maps $L^{[i]}$ are called the \emph{layers} of the neural network. We say that the network has $l$ layers with width $w$, where layers $1,\dotsc,l-1$ are called \emph{hidden layers}.
	
	The set of weights and biases is called the \emph{parameters} of the neural network, given by $\theta \coloneqq \{W^{[i]},b^{[i]}\}_{i = 0,\dotsc,l}$.
	
	Then, the \emph{neural network} associated with the parameters $\theta$ is a map $\soln_\theta \in C^\infty(\symbb{R}^d,\symbb{R})$ given by
	\begin{equation}
	\soln_\theta = L^{[l]}\circ \activation \circ L^{[l-1]} \circ \activation \dotso \circ \activation \circ L^{[1]} \circ \activation \circ L^{[0]},
	\end{equation}
	where we understand the composition with $\activation$ to mean the element-wise application to the layer output.
\end{definition}		 
We limit our analysis to FCFF networks to focus on the effect of the choice of loss function, however our arguments extend to different architectures as well. Further, we limit ourselves to smooth activation functions. It should be noted that a popular choice of activation is the $\operatorname{ReLU}(\cdot)  \coloneqq \max(0,\cdot)$ function, which is not smooth. Our analysis also applies to $\operatorname{ReLU}$ networks.

PINN-based methods are based on minimizing the PDE residuals measured in some discrete approximation of an $L^p$ norm, most commonly for $p=2$.

For $v \in C^1(\tDom\times \domain)$, we define the pointwise PDE residual as 
\begin{equation}
\residual[v] \coloneqq v_t + f(v)_x \in C^0(\tDom \times \domain).
\end{equation}
Let $\soln_\theta$ be a neural network with parameters $\theta$.	Then the classical PINN approach defines the interior contribution to the neural network loss function as a Monte-Carlo approximation of the $L^2$-in-space-time norm of $\residual[\soln_\theta]$. This is done by randomly choosing \textit{collocation points} $\{(t_i,x_i)\}_{i=1}^{N_c}$ with $t_i \in \tDom$, $x_i \in \domain$ and $N_c$ the number of collocation points, and then computing
\begin{equation}
\dLoss_{\mathrm{int}}(\soln_\theta) \coloneqq \frac{1}{N_c} \sum_{i=1}^{N_c} \abs{\residual[\soln_\theta](t_i,x_i)}^2.
\end{equation}
Initial- and boundary data are imposed using further penalty terms, given by
\begin{align}
\dLoss_{\mathrm{ic}}(\soln_\theta) = \norm{\soln_\theta(0,\cdot) - \soln_0}_{2,\mathrm{MC}}^2  &&\text{and}&& \dLoss_{\mathrm{bc}}(\soln_\theta) = \norm{\soln_\theta - g}_{2,\mathrm{MC}}^2,
\end{align}
where $\norm{\cdot}_{2,\mathrm{MC}}$ denotes a Monte-Carlo approximation to the $L^2$-norm on $\domain$ for the initial data and on $(0,T)\times\partial\domain$ for the boundary data. The overall loss is then the weighted sum of interior, initial and boundary loss, $\dLoss(\soln_\theta) \coloneqq \dLoss_{\mathrm{int}} + \lambda \left(\dLoss_{\mathrm{ic}} + \dLoss_{\mathrm{bc}}\right)$ with some user-chosen hyperparameter $\lambda>0$ setting the relative importance of the two terms \cite{WANG_NTK}. As empirically observed in  \cite{Lye_2020}, this approach based on minimizing the $L^2$ norm of the interior residual does not work for learning discontinuous solutions to nonlinear hyperbolic conservation laws.

\subsection{Choice of Norm}
In the following, we present an illustrative argument that shows why classical PINNs cannot learn  such discontinuous solutions based on the example of standing shock solutions to the inviscid Burgers equation. The intuition we develop here easily extends to other nonlinear hyperbolic  conservation laws and moving shocks and highlights a fundamental issue. While the computations are elementary, we could not find an argument of this type in the literature. Instead, the failure of PINNs was treated as an experimental observation so far.

Note that the failure of PINNs is a nonlinear effect and will not occur in the linear case. In fact, in the linear case we expect minimizing $L^p$ norms of the PDE residual  for $1 \leq p \leq \infty$ to work well. For linear advection with constant velocity, by learning a simple linear coordinate transformation, it is clear that the residual can be made to vanish exactly. It remains to learn initial- and boundary data. However, neural networks can easily learn approximations of discontinuous initial data, giving good approximate solutions for linear problems.

Consider that by construction, for continuous activation functions $\activation$, neural networks are continuous as well, so they cannot represent discontinuities exactly. As such, one expects a reasonable neural network approximation to a discontinuous solution to learn smoothed-out versions of shocks, such that it will be close to the exact solution in a pointwise sense outside of a small region around discontinuities and in an $L^2$-sense for the entire domain.

As a concrete example we consider the inviscid Burgers equation, that is the scalar conservation law for the flux function $f(u) = \frac{1}{2}u^2$ and the initial condition 
\begin{equation}
\soln_0(x) = 
\begin{cases}
\begin{aligned}
&1&& \text{for } x<0, \\
-&1&& \text{for } x\geq 0,
\end{aligned}
\end{cases}
\label{eqn:shockSoln}
\end{equation}
having a standing shock as its entropy solution.

As a first scenario we consider time-independent approximations to the solution of this problem, as the exact  solution is time-independent as well and this avoids technical complications.
We consider approximate solutions $\tilde\soln$ which are smooth and satisfy
\begin{equation}
\begin{cases}
\begin{aligned}
&\tilde\soln(x) \in [1-\epsilon,1+\epsilon]&& \text{for } x<-\epsilon, \\
&\tilde\soln(x) \in [-1-\epsilon,-1+\epsilon]&& \text{for } x> \epsilon,\\
&\tilde\soln(x) \in [-1 - \epsilon, 1 + \epsilon] && \text{for } -\epsilon \leq x \leq \epsilon
\end{aligned}
\end{cases}
\label{eqn:shockApprox}
\end{equation}
for some $\epsilon > 0$. For $\epsilon \rightarrow 0$ one easily checks that such approximations converge towards the exact solution in $L^2(\tDom\times\domain)$ with a rate of $\sqrt{\epsilon}$. Our key argument is that these reasonable approximations will, in fact, not produce a small $L^2$ norm of the PDE residual, such that PINNs would not learn such approximations during the optimization process. We show that the $L^2$ norm of the residual scales as $\frac{1}{\sqrt{\epsilon}}$ in this case, i.e. better shock resolution leads to residuals with larger $L^2$ norm.

We consider only the contribution to the residual from $\tDom \times (-\epsilon,\epsilon)$, as the shock is what causes the $L^2$ norm of the residual to grow large and \begin{equation*} \norm{\residual[\tilde\soln]}_{L^2(\tDom\times\domain)} \geq \norm{\residual[\tilde\soln]}_{L^2(\tDom\times(-\epsilon,\epsilon))}.\end{equation*} 

Additionally, Hölder's inequality implies that $\norm{\residual[\tilde\soln]}_{L^1(\tDom\times(-\epsilon,\epsilon))} \leq \sqrt{2\epsilon T}\, \norm{\residual[\tilde\soln]}_{L^2(\tDom\times(-\epsilon,\epsilon))}$. Further, we note that for $\epsilon < 1$, because $\tilde\soln(-\epsilon) > 0$ and $\tilde\soln(\epsilon) < 0$ and $\tilde\soln$ is continuous, there is at least one zero $\bar x$ of $\tilde\soln$ in $(-\epsilon,\epsilon)$. Hence, we have $f(\tilde\soln(\bar x)) = 0$. This behavior close to the discontinuity is what makes the $L^2$ norm of the residual large, because contrary to the true solution which satisfies $f(\soln) = \frac 12$ everywhere, there is a region where $\abs{f(\tilde\soln)_x}$ is large. Then we compute: 
\begin{equation}
\begin{aligned}
\norm{\residual[\tilde\soln]}_{L^1(\tDom\times(-\epsilon,\epsilon))} &= \tInt \int_{-\epsilon}^{\epsilon} \abs{f(\tilde\soln)_x}\dd{x}\dd{t}\\
&\geq \tInt \abs{\int_{-\epsilon}^{\bar x} f(\tilde\soln)_x\dd{x}}\dd{t} +  \abs{\int_{\bar x}^{\epsilon} f(\tilde\soln)_x\dd{x}}\dd{t}\\
&= \tInt \frac{1}{2}\tilde\soln(-\epsilon)^2 + \frac{1}{2}\tilde\soln(\epsilon)^2\dd{t}
\geq T\left(1 - \epsilon\right)^2,
\end{aligned}
\end{equation}
where in the last line we inserted the flux function for the Burgers equation explicitly.

Due to our preliminary considerations the $L^2$ norm of the residual may be estimated by:
\begin{equation}
\label{eqn:badScaling}
\norm{\residual[\tilde\soln]}_{L^2(\tDom\times\domain)}\geq
\norm{\residual[\tilde\soln]}_{L^2(\tDom\times(-\epsilon,\epsilon))}\geq \frac{\norm{\residual[\tilde\soln]}_{L^1(\tDom\times(-\epsilon,\epsilon)}}{2\sqrt{\epsilon T}} \geq \frac{\sqrt{T}(1-\epsilon)^2}{2\sqrt{\epsilon}}.
\end{equation}	
As this estimate shows, for time-independent reasonable approximations to shock solutions of the Burgers equation, reducing $\epsilon$, which intuitively makes the approximation better, increases the $L^2$-norm of its residual.

We are not able to prove a similar estimate for the general case of time-dependent approximations satisfying \eqref{eqn:shockApprox} for each $t$. However, we consider two examples of typical time dependence, which cover a wide range of possible approximations and show that these do not produce a small $L^2$-norm of the residual either. The more technical example is located in Appendix \ref{apx:time_dependent}.

Now we consider approximations whose solution profile is shifted with some time-dependent velocity, but whose shape does not change. This means we consider time-dependent approximations of the form $\tilde\soln(x - s(t))$  with $\tilde\soln$ as in \eqref{eqn:shockApprox} for a smooth function $s$ denoting the overall shift of the solution. Then the time derivative of such approximations is given by $\tilde\soln_t = -s'\tilde\soln_x$. The velocity of the profile is given by $v \coloneqq s'$.
As before, we compute
\begin{equation}
\begin{aligned}
\tInt \int_{-\epsilon + s}^{\epsilon+ s} \abs{f(\tilde\soln)_x +  \tilde \soln_t}\dd{x}\dd{t}
&\geq \tInt \abs{\int_{-\epsilon+s}^{\bar x+s} f(\tilde\soln)_x  - v\tilde\soln_x\dd{x}}\dd{t} +  \tInt \abs{\int_{\bar x+s}^{\epsilon+s} f(\tilde\soln)_x - v\tilde\soln_x\dd{x}}\dd{t}\\
&= \tInt \abs{\underbrace{f\left(\tilde\soln(-\epsilon)\right)-v(t) \tilde\soln(-\epsilon)}_{ \eqqcolon A}} + \abs{\underbrace{f\left(\tilde\soln(\epsilon)\right)- v(t) \tilde\soln(\epsilon)}_{ \eqqcolon B}}\dd{t}\\
& \geq \frac{T}{2}\left(1 - \epsilon\right)^2 + (1 - \epsilon)\norm{v}_{L^1(\tDom)} .
\end{aligned}
\end{equation}
For the last inequality, we use that $\tilde\soln(\epsilon)$ and $\tilde\soln(-\epsilon)$ have opposing sign. Then,  regardless of the sign of $v(t)$, either $A$ or $B$ is the sum of two positive contributions, while the other has mixed sign. For the sum of two positive terms, we may leave out the absolute value and we drop the contribution from the mixed term. The same arguments as before show that the $L^2$ norm of the residual grows as $\frac{1}{\sqrt{\epsilon}}$ again.

We see that for fairly general reasonable approximations, the $L^2$-in-space-time norm of the PDE residual grows large for small $\epsilon$, while the magnitude of the $L^1$ norm does not seem to be connected to $\epsilon$. This shows that both these norms are unsuitable choices for the loss function of PINNs then, because training networks by minimizing the resulting loss functions will not learn approximations with small $\epsilon$.

Instead, one should consider a weaker norm of the residual. We argue that using the {$L^2(0,T;W^{-1,p}(\domain))$- norm ($L^2$-$W^{-1,p}$-norm for short)  with $1 < p < \infty$} is a good choice. Consider again the toy problem as outlined in \eqref{eqn:shockSoln} with a smooth approximation $\tilde\soln$ as in \eqref{eqn:shockApprox}. For technical simplicity we will limit ourselves to approximations that are constant in time, that is $\tilde\soln(t,x) = \tilde\soln(0,x)$ for all $(t,x)$ in $\tDom \times \domain$.

We show that the $L^2$-$W^{-1,p}$ norm of the residual of approximate solutions of the form \eqref{eqn:shockApprox} goes to zero for $\epsilon \to 0$ provided that the sequence is sufficiently non-oscillatory, i.e. if 
{ $\epsilon^{1+ \frac{1}{p}} \norm{f(\tilde\soln)_x}_{L^\infty(\tDom\times (-\epsilon,\epsilon))}$ and $\norm{f(\tilde\soln)_x}_{L^2\left(\tDom;L^p(D\smallsetminus (-\epsilon,\epsilon))\right)}$} both go to zero, which is a desirable property of an approximation.

Fix $1 < p < \infty$ and let $\frac{1}{p} + \frac{1}{q} = 1$. Let $V \coloneqq W^{1,q}_0(\domain)$ with $\norm{\cdot}_V \coloneqq \norm{\cdot}_{W^{1,q}}$, and let $S \coloneqq \{\varphi \in V \colon \norm{\varphi}_V = 1 \}$. Then by definition of the $W^{-1,p}$ norm, we compute
\begin{gather}
\begin{aligned}
\norm{\residual[\tilde\soln]}_{L^2(\tDom;W^{-1,p}(\domain))}^2 
&=\int_{\,\tDom}\left(\sup_{\varphi \in S}  \int_{\domain} f(\tilde\soln)_x \varphi \dd{x}\right)^2\dd{t}\\
&\leq \int_{\,\tDom} \left(\sup_{\varphi \in S}\int_{(-\epsilon,\epsilon)} f(\tilde\soln)_x \varphi \dd{x}\right)^2\dd{t}
+ \int_{\,\tDom}\left(\sup_{\varphi \in S} \int_{D\smallsetminus (-\epsilon,\epsilon)} f(\tilde\soln)_x \varphi \dd{x}\right)^2\dd{t}
\end{aligned}
\label{eqn:h1Norm}
\end{gather}

To estimate the first term of \eqref{eqn:h1Norm}, we rewrite $\varphi(x) = \varphi(0) + \int_0^x \varphi'(s) \dd{s}$ and use the fact that $\abs{(f(\tilde\soln(t,\epsilon)) - f(\tilde\soln(t,-\epsilon))}\leq 4\epsilon$ for the Burgers equation with approximations $\tilde\soln$ as in \eqref{eqn:shockApprox} to conclude that
\begin{gather}
\begin{aligned}
\tInt \left(\sup_{\varphi \in S}\int_{-\epsilon}^{\epsilon} f(\tilde\soln)_x \varphi \dd{x}\right)^2\dd{t} 
&=\!\begin{aligned}[t]
\tInt\left(\sup_{\varphi \in S}\right.&\biggl( \varphi(0) \left(f(\tilde\soln(t,\epsilon)) - f(\tilde\soln(t,-\epsilon))\right)\\
&+\left.{\vphantom{\sup_{\varphi \in S}}}
\int_{-\epsilon}^{\epsilon} f(\tilde\soln)_x \int_0^x \varphi'(s) \dd{s} \dd{x}\biggr)\right)^2\dd{t}
\end{aligned}\\
&\leq  \tInt\left(\sup_{\varphi \in S} \norm{f(\tilde\soln)_x}_\infty \int_{-\epsilon}^{\epsilon} \int_0^x \abs{\varphi'(s)} \dd{s} \dd{x} +  4 \epsilon \varphi(0)\right)^2\dd{t}\\
&\leq  \tInt\left(\sup_{\varphi \in S} \norm{f(\tilde\soln)_x}_\infty \int_{-\epsilon}^{\epsilon} {\abs{x}}^{\frac{1}{p}}\norm{\varphi'}_{L^q(\domain)}\dd{x} +  4 \epsilon \varphi(0)\right)^2\dd{t}\\
&=  T \sup_{\varphi \in S}\left(\norm{f(\tilde\soln)_x}_\infty\left(\frac{2}{1 + 1/p}\epsilon^{1 + \frac{1}{p}} \right) \norm{\varphi'}_{L^q(\domain)} +  4 \epsilon \varphi(0)\right)^2\\
&\leq   T \left(\norm{f(\tilde\soln)_x}_\infty\left(\frac{2}{1 + 1/p}\epsilon^{1 + \frac{1}{p}} \right) +  4 \epsilon c_S\right)^2.
\end{aligned}
\end{gather}
The last inequality uses that, by definition, $\norm{\varphi' }_{L^q(\domain)} \leq 1$ and that in one dimension $W^{1,q}$ with $q>1$ is continuously embedded in $L^\infty$, with embedding constant denoted by $c_S$.

Estimating the second term of \eqref{eqn:h1Norm} is straightforward. Note that by definition of $\varphi$ and applying Hölder's inequality, one sees that 
\begin{equation}
\int_{\,\tDom}\left(\sup_{\varphi \in S} \int_{\domain\smallsetminus (-\epsilon,\epsilon)} f(\tilde\soln)_x \varphi \dd{x}\right)^2\dd{t} \leq \int_{\,\tDom}\left( \int_{\domain\smallsetminus (-\epsilon,\epsilon)} f(\tilde\soln)^p_x  \dd{x}\right)^{\frac{2}{p}}\dd{t} =  \norm{f(\tilde\soln)_x}_{L^2\left(\tDom;L^p(D\smallsetminus (-\epsilon,\epsilon))\right)}^2.
\end{equation}
Combining the estimates for each of the two previous terms, we may estimate the residual in this norm by
\begin{equation}
\label{eqn:h1Estimate}
\norm{\residual[\tilde\soln]}_{L^2(\tDom;W^{-1,p}(\domain))}^2 \leq T \left(\norm{f(\tilde\soln)_x}_\infty\left(\frac{2}{1 + 1/p}\epsilon^{1 + \frac{1}{p}} \right) +  4 \epsilon c_S\right)^2 + \norm{f(\tilde\soln)_x}_{L^2\left(\tDom;L^p(D\smallsetminus (-\epsilon,\epsilon))\right)}^2.
\end{equation}
The last term of estimate \eqref{eqn:h1Estimate} is not connected to the resolution of the shock.  Comparing the two terms inside the square, because $\norm{f(\tilde\soln)_x}_\infty \geq \frac{1}{2\epsilon}$, the first term scales with $\epsilon^{\frac{1}{p}}$ while the second term scales with $\epsilon$, so the first term is the dominant contribution for small $\epsilon$. Larger values of $p$  produce steeper gradients of the loss with respect to $\epsilon$ when $\epsilon$ is small and will incentivize finer shock resolution. However, choosing larger $p$ may not always be better, because larger $p$ have smaller gradients for large values of $\epsilon$, which occur at the beginning of training neural networks.

Overall, the estimate shows that the $L^2$-$W^{-1,p}$ norm with $1 < p < \infty$ is a good choice of norm for the PDE residual in the neural network loss function, because 
good approximations of the exact solution correspond to small losses.
Approximations oscillating in $\tDom\times (-\epsilon,\epsilon)$ will lead to larger values of $\norm{f(\tilde\soln)_x}_\infty$. Such approximations do not guarantee a small norm of the residual, showing that the goal functional discourages spurious oscillations.

\begin{remark}
	The above arguments were outlined in the case of a stationary shock, but can be adapted to cover moving shocks. The estimates employed are done pointwise for each $t \in \tDom$, such that the only difference is that the interval $(-\epsilon,\epsilon)$ now moves with the proper shock velocity instead of being fixed.
\end{remark}

The original wPINNs fit into the previous analysis. They introduce a \textit{(Kruzhkov) entropy residual} \cite{wPINNS}, based on enforcing entropy inequalities, which is then turned into a loss function by measuring it in something similar to the $L^2(0,T;W^{-1,2}(\domain))$ norm. For $\bar S \coloneqq \{\testSoln \in C_c^\infty(\tDom \times \domain) : \testSoln\geq 0, \norm{\partial_x \testSoln}_{L^2(\tDom\times\domain)} = 1\}$, the loss function on the continuous level is given by
\begin{equation}
\Loss_{\mathrm{ent}}^K(\soln_\theta) = \sup_{\testSoln\in \bar S} \max_{c \in \symbb{R}}\int_{\,\tDom}\int_{\domain} \testSoln\partial_t\abs{\soln_\theta - c} - \sgn(\soln_\theta-c)(f(\soln_\theta)-f(c))\testSoln_x \dd{x} \dd{t}.
\end{equation}
Note that this form of loss contains a nested maximization problem over $c \in \symbb{R}$, which is approximated in \cite{wPINNS} by sampling some discrete set of $c$ and computing the maximum directly.

In contrast to our previous arguments, $\Loss_{\mathrm{ent}}^K$ is not based on computing a norm of the PDE residual directly, however the Kruzhkov entropy inequalities already encode the PDE itself. Being based on entropy conditions, this choice of loss ensures that one finds the unique entropy solution.
	\section{Modified wPINNs and Extensions}
\label{sec:weakNorms}
\subsection{Weak Norm Estimation}
Our previous arguments show that classical PINNs cannot succeed for hyperbolic conservation laws and it is instead required to use PINNs based on dual norms for this problem class.
We use a novel approach for computing dual norms using neural networks by solving a dual elliptic problem. This improves the computational efficiency of weak PINNs by giving better approximations of the dual norm with less training. The modified loss functional accelerates learning in terms of the number of epochs the neural networks require to be trained to comparable accuracy, or achieve higher accuracy after equal number of epochs.

To fix notation, we first describe the loss components on an analytical level. We replace the Kruzhkov entropy formulation from the original wPINNs by the PDE residual and an entropy residual based on a single strictly convex entropy, because systems of hyperbolic conservation laws commonly only have one physically motivated entropy. Note that this is independent of our approach to computing dual norms, which can be applied to the Kruzhkov entropy formulation of the original wPINNs as well.

The interior PDE residual loss measures the PDE residual in the $L^2$-$W^{-1,p}$ norm. Recall that $S = \{\varphi \in V \colon \norm{\varphi}_V = 1 \}$ with $V = W^{1,q}_0(\domain)$ and $q$ is the Hölder-conjugate exponent of $p$. The interior PDE residual loss is then given by 
\begin{equation}
\Loss^u_{\mathrm{int}}(\soln_\theta) = \int_{\,\tDom}\biggl(\sup_{\testSoln \in S} \int_{\domain} \left(\partial_t \soln_\theta + \partial_x f(\soln_\theta)\right) \testSoln\dd{x}\biggr)^2\dd{t}.
\end{equation} 
This loss is minimized by any weak solution, even ones which are not entropy solutions, so enforcing an entropy inequality is required. In \cite{wPINNS}, a numerical example approximating a non-entropic weak solution is shown.  Provided $f$ is strictly convex, a single strictly convex entropy-entropy flux pair is sufficient to enforce entropy admissibility, instead of the family of Kruzhkov entropies. Thus, we fix one such entropy-entropy flux pair $(\entropy,\entflux)$ and introduce the \textit{entropy loss}, measured in the $L^2$-$W^{-1,p}$ norm, as 
\begin{equation}
\Loss_{\mathrm{ent}}^\eta(\soln_\theta) = \int_{\,\tDom} \biggl(\sup_{\testEnt \in S} \int_{\domain} \left(	\entropy(\soln_\theta)_t + \entflux(\soln_\theta)_x\right)^{\oplus}\testEnt \dd{x}\biggr)^2\dd{t},
\end{equation} 
where we take the positive part denoted by $(\cdot)^\oplus$ since the entropy condition is an inequality. This approach avoids the nested maximization over $c\in \symbb{R}$ of the original wPINN approach in favor of a second, independent, maximization problem with another adversarial neural network that one trains using gradient ascent. 

We now discuss our strategy for approximating the $W^{-1,p}$ norm, whose definition takes a supremum over all $\varphi \in S$.
When parameterizing approximations to $\varphi$ using a neural network, it is unclear how to enforce the normalization of the test function through the network design, so one has to compute the $W^{1,q}$ norm (using Monte-Carlo sampling) and divide by this. We believe that the division by such a norm makes the maximization problem more difficult. Also, by construction, this formulation is invariant under scalar rescaling of the network approximating the dual norm. As the updates during gradient ascent change the normalization of the network, this undetermined degree of freedom may reduce the speed and accuracy of the norm estimation.

To remedy these downsides, we propose computing the $W^{-1,p}$ norm by learning the solution to a $q$-Laplace problem with $\frac 1p + \frac 1q = 1$. We  identify the dual space $W^{-1,p}(\domain) = \left(W^{1,q}_0(\domain)\right)'$. Given a functional $v \in W^{-1,p}(\domain)$ for $1 < p < \infty$, we see that determining
\begin{equation}
w \in W^{1,q}_0(\domain) \text{ with } \int_\domain \abs{\grad w}^{q-2}\grad w \grad \zeta \dd{x} = \int_\domain v \zeta \dd{x} \quad \text{for all } \zeta \in W^{1,q}_0(\domain),
\end{equation}
one has $\norm{v}^p_{W^{-1,p}(\domain)} = \abs{w}^q_{W^{1,q}(\domain)} \coloneqq \norm{\grad{w}}^q_{L^q(\domain)}$. For a derivation of the norm relationships between $v$ and $w$, see \cite{Dinca_pLap}. Note that due to the zero-boundary conditions, by Poincaré inequality the $W^{1,q}$-seminorm is a proper norm on $W^{1,q}_0(\domain)$. The solution of this dual problem may also be characterized as the maximizer of the energy functional 
\begin{equation}
I(w) = \int_\domain v w \dd{x} - \frac{1}{q}\int_{\domain} \abs{\grad{w}}^q \dd{x}.
\end{equation}
We use this characterization to train our test functions. In the case $q=2$ maximizing $I(w)$ on the set of neural functions is known as the “Deep Ritz Method” \cite{E2018} and is robust and effective even for low-regularity right-hand-sides $v$ for the Poisson problem \cite{Chen_2020}.

\subsection{Loss Definition}
Next, we discretize the continuous loss functional using Monte-Carlo sampling. We introduce
interior collocation points $\collocPoints_{\mathrm{int}} \coloneqq \{(t_i,x_i)\}_{i=1}^{N_{\mathrm{int}}}$ with $t_i \in \tDom$ and $x_i \in \domain$ and $N_{\mathrm{int}}$ the number of collocation points. Similarly, we have  $S_\mathrm{ic}$ sampling $\domain$ for the initial data and $S_\mathrm{bc}$ sampling $(0,T)\times\partial\domain$ for the boundary data.

We work with three neural networks $\soln_\theta$, $\tilde\testSoln_\chi $ and $\tilde\testEnt_\nu$ with weights denoted by $\theta, \chi \text{ and } \nu$. The latter two networks approximate the solutions to the dual problems as above, so one needs to enforce zero boundary conditions through either soft- or hard constraints \cite{Chen_2020}. Soft constraints add an additional loss term penalizing the violation of boundary conditions. Hard constraints encode the boundary condition exactly into the network design. We opt for hard constraints, because this speeds up learning and improves accuracy of PINNs  \cite{Chen_2020}. We enforce the boundary conditions through a cutoff function $w:\domain \rightarrow \symbb{R}^+$ with $w\rvert_{\partial\domain} = 0$, $w' \rvert_{\partial\domain} \neq 0$ and $w(x) >0$ for all $x$ in the interior of $\domain$. Then the functions $\testSoln_\chi(x,t) \coloneqq \tilde\testSoln_\chi(x,t)w(x)$ and $\testEnt_\nu(x,t) \coloneqq \tilde\testEnt_\nu(x,t)w(x)$ satisfy zero boundary conditions.

To simplify notation we define
\begin{equation}
\resfun_{\mathrm{PDE}}(\soln_\theta;\testSoln_\chi) := (\partial_t\soln_{\theta})\testSoln_\chi - f(\soln_\theta) (\partial_x\testSoln_{\chi}).
\end{equation}
This definition corresponds to the integrand of the PDE residual loss after integration by parts in space, as the boundary contributions vanish because the network $\testSoln_\chi$ has zero boundary conditions in space. The integration by parts when defining the loss for numerical simulations is also done in \cite{wPINNS} and gives good results, so we adopt this approach. The overall PDE residual loss is then discretized as
\begin{equation}
\dLoss_{\mathrm{PDE}}^\soln \coloneqq \frac{1}{N_{\mathrm{int}}}\biggl(\sum_{\collocPoints_{\mathrm{int}}} \resfun_{\mathrm{PDE}}(\soln_\theta;\testSoln_\chi) (t_i,x_i)- \frac{1}{q} \abs{\grad_x\testSoln_\chi}^q(t_i,x_i)\biggr).\label{eqn:mixedRes}
\end{equation}
We solve the elliptic problem in space using a single network with the time $t$ as an additional input parameter, such that it learns the solution for the dual problem simultaneously on the entire time interval. This is equivalent to the continuous definition of the loss, but more effective on the discrete level.

Note that, when $\testSoln_\chi$ (approximately) solves the aforementioned dual elliptic problem with the PDE residual as right-hand side, \eqref{eqn:mixedRes} corresponds to the Monte-Carlo approximation of  $\int_{[0,T)}\int_\domain \residual[\soln_\theta] \testSoln_\chi - \frac{1}{q} \abs{\grad_x{\testSoln_\chi}}^q \dd{x}\dd{t} = \frac{1}{p} \int_{\,\tDom} \int_\domain \abs{\grad_x \testSoln_\chi}^q\dd{x}\dd{t}$ . We compute the approximation of the $L^2$-$W^{-1,p}$ norm using equation \eqref{eqn:mixedRes} because the dependence on the parameters $\theta$ of the network $\soln_\theta$ is explicit then, as required for automatic differentiation.

Further, we have the \textit{entropy residual} for a given entropy-entropy flux pair $(\entropy,\entflux)$ defined as 
\begin{equation}
\resfun_{\mathrm{ent}}(\soln_\theta;\testEnt_\nu;(\eta,q)) \coloneqq \left(\entropy(\soln_\theta)_t + \entflux(\soln_\theta)_x\right)^{\oplus}\testEnt_\nu,
\end{equation}
which we measure again in the squared $L^2$-$W^{-1,p}$ norm as  
\begin{equation}
\dLoss_{\mathrm{ent}}^\eta \coloneqq \frac{1}{N_{\mathrm{int}}}\biggl(\sum_{\collocPoints_{\mathrm{int}}}\resfun_{\mathrm{ent}}(\soln_\theta;\testEnt_\nu) (t_i,x_i) - \frac{1}{q} \abs{\grad_x\testEnt_\nu}^q(t_i,x_i)\biggr).\label{eqn:entRes}
\end{equation}
Putting these parts together we consider the following interior loss function: 
\begin{align}
\dLoss_{\mathrm{int}} = 	\dLoss_{\mathrm{ent}}^\eta + \dLoss_{\mathrm{PDE}}^\soln
\end{align}
Note that the interior loss consists of two additive contributions, the PDE residual loss \eqref{eqn:mixedRes} and the entropy loss \eqref{eqn:entRes}, each of which again consist of two parts. The PDE residual loss $	\dLoss_{\mathrm{PDE}}^\soln$ has no dependence on the adversarial network $\testEnt_\nu$, and the spatial gradient part depends only on the network $\testSoln_\chi$, while the entropy loss $	\dLoss_{\mathrm{ent}}^\eta$ does not depend on $\testSoln_{\chi}$ and its spatial gradient part depends only on $\testEnt_\nu$.

For the overall loss we add in initial- and boundary contributions, for Dirichlet boundary data $g$ on the spatial boundary of the domain and initial data $u_0$ at $t=0$. This gives the overall loss function
\begin{equation}
\begin{aligned}
\dLoss(\soln_\theta,\testSoln_\chi,\testEnt_\nu) =& \dLoss_{\mathrm{int}}(\soln_\theta,\testSoln_\chi,\testEnt_\nu) \\&+ \lambda\left(\frac{1}{N_\mathrm{ic}} \sum_{\collocPoints_{\mathrm{ic}}} (\soln_\theta(0,x_i) - \soln_0(x_i))^2 + \frac{1}{N_\mathrm{bc}} \sum_{\collocPoints_{\mathrm{bc}}} (\soln_\theta(t_i,x_i) - g(t_i,x_i))^2\right),
\end{aligned}
\label{eqn:Loss}
\end{equation}
with parameter $\lambda > 0$.
These additional terms depend only on the solution network $\soln_\theta$. The decomposition of the loss into additive terms independent of some of the networks allows avoiding some computations to improve efficiency. For example, computing a gradient of the loss with respect to only the parameters $\theta$ does not require computation of the spatial gradients  of the networks $\testSoln_{\chi}$ and $\testEnt_\nu$ in \eqref{eqn:mixedRes} and \eqref{eqn:entRes}. The parameter $\lambda$ may be chosen to balance gradients of the interior and boundary contributions during training \cite{wang_gradientflow}.

\subsection{Extension to Weak Boundary Conditions}
\label{sec:weakBC_theo}
In this section, we show how to incorporate Dirichlet boundary conditions in a weak sense for scalar conservation laws using the framework from \cite{kondo2001measure}. Let $u$ be a weak solution of \eqref{eqn:cLaw}--\eqref{eqn:iCond} in the sense of \eqref{eqn:wcLaw}.

Then, the bounded boundary data $u(t,x) = g(t,x)$ on $(0,T) \times \partial\domain$  need to be understood in a weak sense, i.e. the value of $g$ is only attained, if it lies on the inflow part of the boundary, while on the outflow part of the boundary, the boundary data need not be attained. 

Mathematically, for a fixed convex  entropy-entropy-flux pair $(\eta,q)$ we express this as
\begin{equation}
\int_{\tDom} \int_{\partial\domain} \Bigl( \left( q(g) - q(u)\right)\cdot n + \eta'(g)\left(f(u)-f(g)\right)\cdot n \Bigr) \varphi \dd{x} \dd{t} \leq 0,
\label{eqn:wBound}
\end{equation}
for all smooth $\varphi \geq 0$, where $n$ represents the outer unit normal on $\partial\domain$. We write the integral expression from \cite{kondo2001measure} despite working in one space dimension to make the underlying structure more clear.

To enforce \eqref{eqn:wBound}, we  replace the standard PINNs boundary term $\norm{u_\theta - g}_{L^2((0,T)\times\partial\domain)}^2$ by
\begin{equation}
\Loss_{\mathrm{bc}}(u_\theta) = \norm{\Bigl( \left( q(g) - q(u_\theta)\right)\cdot n + \eta'(g)\left(f(u_\theta)-f(g)\right)\cdot n \Bigr)^\oplus}_{L^2((0,T)\times\partial\domain)}^2
\end{equation}
The discretization of this term proceeds as previously, by approximating the $L^2$-norm using Monte-Carlo integration on $(0,T)\times\partial\domain$.

\subsection{Causality}
Solutions to hyperbolic conservation laws have causal structure, meaning that information only propagates forward in time, with a limited speed. In contrast, PINNs do not enforce this structure.
During numerical experiments for weak boundary conditions, we observed shocks being moved outside of the computational domain at small values of $t$ instead of propagating at the correct speed.

We hypothesize that this is related to the lack of causality of PINNs. Shocks contribute the majority of the neural network loss, while the magnitude of the residual is  significantly smaller away from shocks. On longer time domains, because gradient descent locally minimizes the loss in the parameters of the network, incurring a large residual error at early times $t$ by removing shocks from the domain and having a smooth approximate solution on the rest of the domain may appear a promising descent direction. However, this leads to completely incorrect approximations due to getting stuck in a local loss minimum. This is due to PINNs not enforcing causality. If causality were enforced, making a large error at early times would no longer be “encouraged” through a potential reduction in loss at later times.

Note that this problem may occur with strong boundary data enforced using a penalty as well, when the weight of the boundary loss is too small. Increasing the weight makes violating the boundary conditions prohibitively expensive and can somewhat mitigate this behavior. However, the option of committing a large error at early times and then learning some unrelated smooth solution that is compatible with the boundary data still remains.  With weak boundary conditions, changing the weighting of the loss contributions is ineffective. Weak boundary conditions do not penalize the network for moving the shock from the domain, so increasing their weight does not change the situation. Only the residual loss penalizes the network for making a large error at early times, and weighting the residual more also increases projected payoff of having a smooth solution after making a large error.

To enforce some level of causality, we subdivide the time domain into smaller slices and subsequently train neural networks on each time slice. Details can be found in Section \ref{sec:weakBC_num} where we apply our time-slicing approach. This approach enforces causality between each of the time slices, while on each time slice there is no causality.
Another possible approach would be “causal training” \cite{wang2022causality}, that continuously discounts the loss at later times $t$ based on the size of loss at earlier times, discouraging large errors for small $t$. However, we did  not try this alternative.
	\section{Numerical Results - Performance Comparison}
\label{sec:numResults}
In this section we present several numerical experiments to compare the performance of the modified wPINNs against the original wPINNs. We train neural networks for initial data corresponding to a standing shock, a moving shock, a rarefaction wave from Riemann problem initial data and sine initial data, to match the range of examples covered in \cite{wPINNS}.
\subsection{Algorithm and Implementation Details}%
\label{sec:implementation}%
First, we describe our training procedure which we use for the performance comparison with the original wPINNs. Further notes and explanations detailing the implementation are located in Appendix \ref{apx:implementation} to keep the paper self-contained and ensure reproducibility of our numerical tests. Our python implementation is located at \url{https://git-ce.rwth-aachen.de/aidan.chaumet/wpinns}.
\begin{algorithm}[ht]
	\caption{Neural Network Training algorithm}
	\label{alg:GAN}
	\begin{algorithmic}[1] 
		\Require{Initial data $u_0$, boundary data $g$, flux function $f$, strictly convex entropy-entropy flux pair $(\eta,q)$, Hyperparamters: $\tau_{\mathrm{min}},\tau_{\mathrm{max}},N_\mathrm{max},\lambda,\gamma$}
		\Ensure{Best networks $\soln_\theta^b,\testSoln_\chi^b,\testEnt_\nu^b$}
		\State{Initialize the networks $\soln_\theta,\testSoln_\chi,\testEnt_\nu$}			
		\State{Initialize performance metric $\dLoss_{\mathrm{avg}} \gets \dLoss(\soln_\theta,\testSoln_\chi,\testEnt_\nu)$}
		\State{Initialize best performance metric $\dLoss_{\mathrm{best}} \gets \infty$}
		\State{Generate collocation points $\collocPoints_{\mathrm{int}},\collocPoints_{\mathrm{ic}},\collocPoints_{\mathrm{bc}}$}
		\For{$ep = 1,...,N_\mathrm{ep}$}
		\For{$k = 1,...,N_\mathrm{max}$}
		\State{Compute $\dLoss(\soln_\theta,\testSoln_\chi,\testEnt_\nu)$} \label{alg:loss1}
		\State{Update $\chi \gets \chi + \tau_{\mathrm{max}}\grad_\chi{\dLoss(\soln_\theta,\testSoln_\chi,\testEnt_\nu)}$ }
		\State{Update $\nu \gets \nu + \tau_{\mathrm{max}}\grad_\nu{\dLoss(\soln_\theta,\testSoln_\chi,\testEnt_\nu)}$ }
		\EndFor
		\State{Compute $\dLoss(\soln_\theta,\testSoln_\chi,\testEnt_\nu)$} \label{alg:loss2}
		\State{Update $\theta \gets \theta - \tau_{\mathrm{min}}\grad_\theta{\dLoss(\soln_\theta,\testSoln_\chi,\testEnt_\nu)}$ }
		\State{Update performance indicator $\dLoss_{\mathrm{avg}} \gets (1 - \gamma) \dLoss_{\mathrm{avg}} + \gamma \dLoss(\soln_\theta,\testSoln_\chi,\testEnt_\nu) $}
		\If{$\dLoss_{\mathrm{avg}} < \dLoss_{\mathrm{best}}$}
		\State{$\dLoss_{\mathrm{best}}\gets \dLoss_{\mathrm{avg}}$}
		\State{Save best networks $\soln_\theta^b,\testSoln_\chi^b,\testEnt_\nu^b \gets \soln_\theta,\testSoln_\chi,\testEnt_\nu$}
		\EndIf
		\EndFor
	\end{algorithmic}
\end{algorithm}

Our training procedure is described in Algorithm \ref{alg:GAN}. It is essentially a standard generative adversarial neural network (GAN) training procedure including best model checkpoints, and as such many of the usual heuristics in the training of such network architectures apply.
Note that lines \ref{alg:loss1} and \ref{alg:loss2} of the algorithm compute the loss function depending on the networks. These steps generate the computational graphs to backpropagate gradients using automatic differentiation for weight updates. We write these steps to represent a practical implementation in machine learning frameworks.
Note that the list of hyperparameters given is not exhaustive.

\subsection{Comparison Setup}  
Now, we describe our methodology for comparing the performance of the original and modified wPINNs. For the entirety of this section we consider $p = q = 2$ in the modified loss because this is most closely related to the original wPINNs.

As a measure of computational cost, we count the number of training epochs. The original wPINNs and the modified wPINNs require similar computational effort per epoch, because the most expensive part of each epoch is evaluating the PDE residual numerically which occurs equally often in original and modified wPINNs. We do not compare wall-clock time, which may differ due to implementation details or choice of machine learning framework.

We make no comparison to classical PINNs because these cannot produce good approximate solutions. We also do not compare to conventional numerical methods because it is obvious that in one or two dimensions neural networks cannot compete, while in high dimensions conventional mesh-based methods are obviously unsuitable, so a like-to-like comparison is not possible.

We assess the accuracy of approximations $\soln_\theta$ using the relative space-time $L^1$ error given by
\begin{equation}
\relerr(\soln_\theta) = \frac{\int_{\,\tDom \times \domain}\abs{\soln_\theta(t,x) - \soln(t,x)}\dd{x}\dd{t}}{\int_{\,\tDom \times \domain}\abs{\soln(t,x)}\dd{x}\dd{t}},
\end{equation}
where $\soln$ denotes the exact solution, in case it is easy to compute. The integrals are approximated using Monte-Carlo integration with $2^{17}$ points in the space-time domain. In cases where the exact solution is not readily available, we approximate it using a high-resolution finite volume scheme from PyClaw \cite{pyclaw-sisc}.

We perform an ensemble training procedure as in \cite{wPINNS}, however we do not investigate as large a range of hyperparameters. Instead, we try several hyperparameter configurations until we find parameters that give good results for the original wPINN algorithm and then use the same hyperparameters for training both approaches. Although there is no direct corresponding network for $\testEnt_\nu$ in the original wPINN approach, we match this network's hyperparameters to those of $\testSoln_\chi$.  This keeps the two training procedures as comparable as possible.

As a learning rate schedule, we reduce the learning rate linearly by half over the total number of epochs. While we think that a more aggressive learning rate schedule would benefit the modified wPINNs, to be able to check our results against the original wPINN results from \cite{wPINNS} we maintain their choice.

For both original and modified wPINNs, we track $\relerr(\soln_\theta)$ for each network of the ensemble individually over the training epochs for comparison, sampled at multiple epochs during training. 

Additionally, we also compare the quality of the \emph{average network} predictions
\begin{equation}
\soln^*(t,x) = \frac{1}{N^*} \sum_{i=1}^{N^*} u_{\theta,i}(t,x),
\end{equation}
for an ensemble  $\{u_{\theta,i}\}_{i=1}^{N^*}$ consisting of $N^*$ neural networks trained for the same hyperparameters, but different randomly chosen initial weights and biases. We compute $\relerr(\soln^*)$ only at the final epoch of training.

Summarizing our numerical results,  we list the average network prediction error in Table \ref{tab:shocksL1}. In particular, the table shows that we reproduce the accuracy achieved by the original wPINNs in \cite{wPINNS} and that the modified wPINNs match or outperform the original wPINNs in all examples, with significant gains especially for the less synthetic examples. In the following sections, we discuss each of the numerical experiments in more detail.
\begin{table}[hb]
	\caption{{Comparison of the r}elative error $\relerr(u^*)$ of the average network prediction at the {different times during} training. The modified wPINNs are two times more accurate than the original wPINNs after less than $30\%$ of the total epochs for the rarefaction wave and over 10 times more accurate for the sine initial data.}
	\label{tab:shocksL1}
	\begin{tabular}{|l|c|c|}
		\hline
		&original & modified \\\hline
		Standing shock  (2000 Epochs)& $0.105\%$&$0.098\%$\\\hline
		Moving shock  (2000 Epochs)& $1.46\%$& $1.36\%$ \\\hline 
		Rarefaction Wave (800 Epochs)& $2.94\%$& $1.57\%$ \\\hline 
		Rarefaction Wave (3000 Epochs)& $1.45\%$& $1.18\%$ \\\hline 
		Sine Wave (10000 Epochs)& $21.4\%$& $1.53\%$ \\\hline 
		Sine Wave (75000 Epochs)& $5.03\%$& $1.19\%$ \\\hline 
	\end{tabular}
\end{table}
\subsection{Standing Shock}
We begin by considering the Burgers equation on $[0,0.5) \times [-1,1]$ for the initial datum
\begin{equation}
\soln_0(x) = \begin{cases} 1 \quad&\text{for } x \leq 0, \\ -1 \quad&\text{for } x>0.\end{cases} \qquad 
\end{equation}
The corresponding exact solution is given by
\begin{equation}
\soln(t,x) = \begin{cases} 1 \quad&\text{for } x \leq 0, \\ -1 \quad&\text{for } x>0,\end{cases} \qquad 
\end{equation}
that is, a standing shock at $x=0$.

We use the following hyperparameters which are the same for both procedures: $\tau_{\mathrm{min}} = 0.01, \tau_{\mathrm{max}}= 0.015$, $N_{\mathrm{max}} = 8$, $ N_\mathrm{ep} = 2000$. The networks $\soln_\theta$ in both algorithms have $l = 6$ layers of width $w = 20$ with $\activation(\cdot) = \tanh(\cdot)$ activation. The networks $\testSoln_{\chi}$ and $\testEnt_\nu$ use $4$ layers of width $10$ with $\tanh$ activation function.
We fix the convex entropy-entropy flux pair $\eta(u) = \frac{1}{2}u^2$ and $q(u) = \frac{1}{3}u^3$ for  all following examples.

For the original wPINNs we adopt the penalty parameter $\lambda = 10$, as fixed during the numerical experiments in \cite{wPINNS}, while for the modified wPINNs $\lambda = 1$ seems to work better. Because we use different loss functions it is not sensible to keep this parameter matched. We set $\gamma = 0.3$ for the exponential averaging of past loss values. The original wPINN algorithm also requires a choice of reset frequency for the adversarial network which we choose as $r_f = 0.05$. Lastly, we compute $N^* = 16$ retrainings.

For both the original and modified wPINNs we use $N_\mathrm{int} = 16384$ uniformly randomly sampled interior collocation points, $N_\mathrm{ic} = 4096$ inital collocation points and $N_\mathrm{bc} = 4096$ boundary collocation points. These points are randomly chosen anew for each retraining of neural networks and choice of algorithm, but stay fixed during the training of each single network.

The average network prediction error is listed in table \ref{tab:shocksL1}. In this example we see that the performance of both methods is approximately the same. We give several reasons why this can be mainly attributed to the specific example below.	
\begin{figure}[b]
	\begin{subfigure}{0.45\textwidth}
		\includegraphics[width=\textwidth]{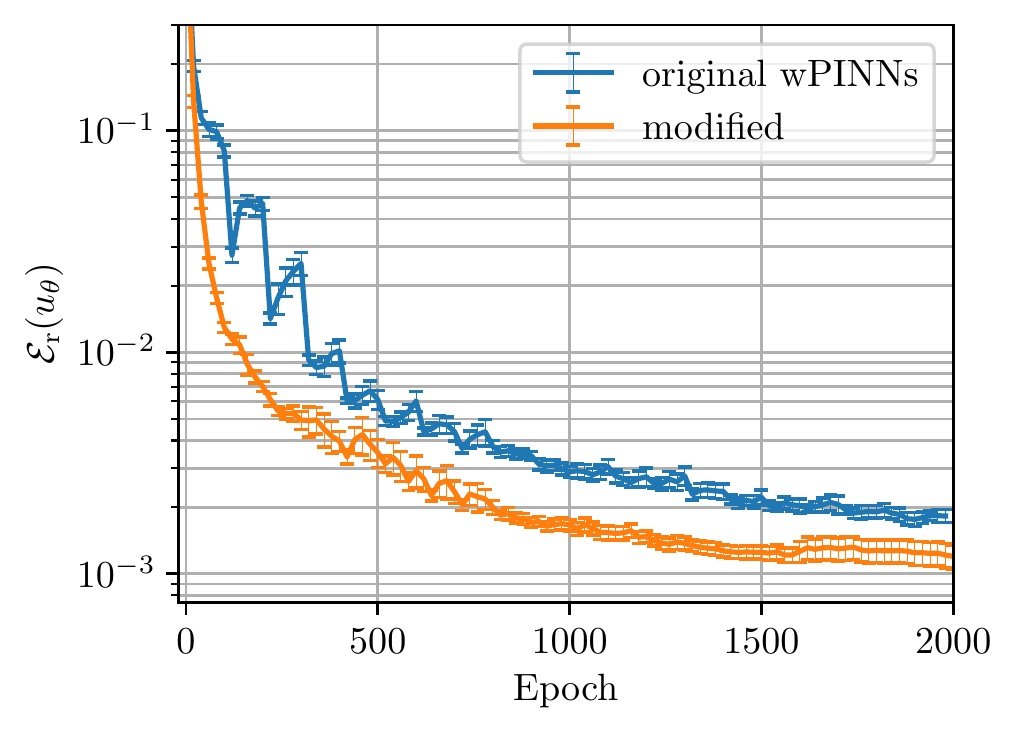}
		\caption{Standing shock}
		\label{fig:standing}
	\end{subfigure}
	\begin{subfigure}{0.45\textwidth}
		\includegraphics[width=\textwidth]{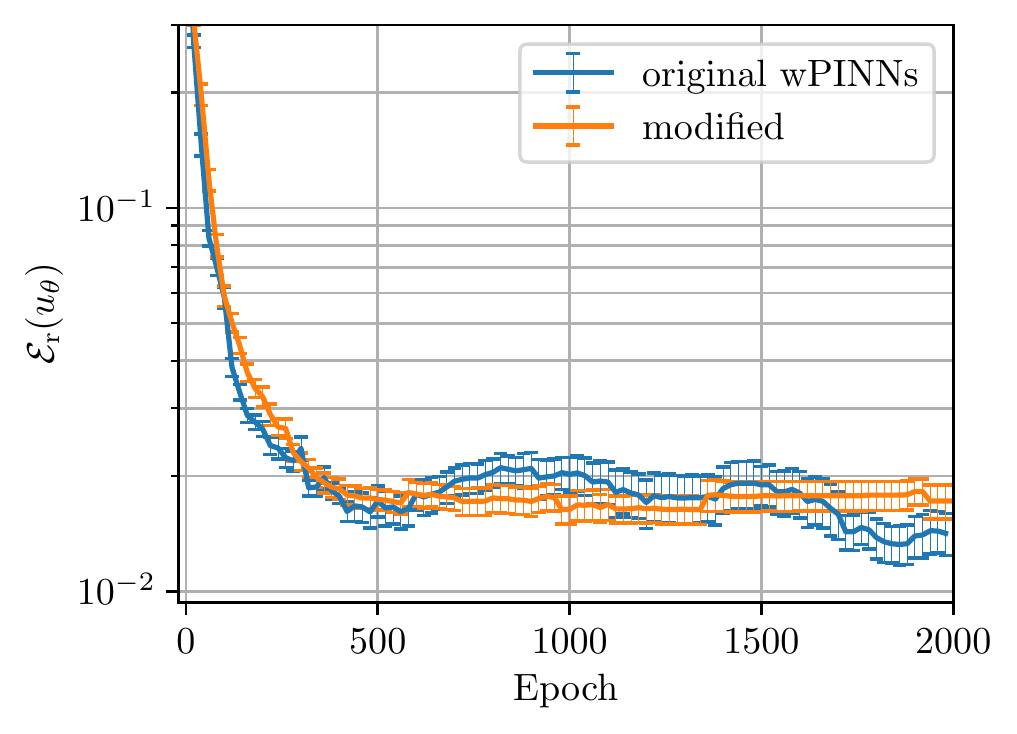}
		\caption{Moving shock}
		\label{fig:moving}
	\end{subfigure}
	\caption{Evolution of the mean value and standard error of $\relerr$ for the respective network ensembles. }
	\label{fig:shocksL1}
\end{figure}

We plot the individual network's errors over the training epochs for the standing shock in Figure \ref{fig:standing}. The error bars show the standard error of the mean, that is $\frac {\sigma}{\sqrt{N^*}}$ with $\sigma$ the standard deviation of the relative $L^1$ errors in the ensemble.	
There are a few small differences in the evolution of $\relerr$, such as the modified wPINNs giving a “smoother” evolution of the relative error over the epochs and having generally smaller statistical variations.

We attribute the smoother evolution of the relative error to the more stable norm estimation of our loss functional by avoiding the division by the $W^{1,2}$-seminorm of the approximating network and possibly the adjusted checkpointing procedure as well, as described in appendix \ref{apx:implementation}.

Confirming our expectations, we see that both original and modified wPINNs perform well, showing that both methods are suitable for approximating shock solutions. The difference between the two methods is small in this example. {The limiting factor for accuracy seems to be the network architecture. As depth and width are identical for either test, the accuracy will thus be similar.} As the initial data is already a shock, the loss contributions from initial and boundary data heavily influence the overall learning of the solution profile. The additional PDE residual contribution only needs to ensure that the shock is placed at the right position, and because the shock is stationary this is straightforward, giving only a small contribution to the total loss.

The initial- and boundary loss functions are identical for original and modified wPINNs. During training one can check that the majority of the overall loss for either approach is due to these contributions for this example. The different choice of penalty parameter $\lambda$ only has a minor influence on the training procedure because the Adam algorithm for minimizing the loss is invariant with respect to scalar rescaling of the entire loss, which consists mainly of the initial- and boundary loss. As such, the original and modified wPINNs are almost equivalent for this example and other situations where the shape of the solution does not differ much compared to the initial data.
\subsection{Moving Shock}
For completeness, we also consider a moving shock example. We choose
\begin{equation}
\soln_0(x) = \begin{cases} 1 \quad&\text{for } x \leq 0, \\ 0 \quad&\text{for } x>0,\end{cases} 
\end{equation}
as an initial datum leading to the solution
\begin{equation}
\soln(t,x) = \begin{cases} 1 \quad&\text{for } x \leq \frac{t}{2}, \\ 0 \quad&\text{for } x>\frac{t}{2},\end{cases}
\end{equation}
which is a moving shock starting from $x=0$, moving to the right with a speed of 0.5.

We choose the same training setup and hyperparameters as for the standing shock.
As expected, the results for this example are structurally the same as in the standing shock example. We plot the results for the moving shock in Figure \ref{fig:moving}, and reference the average network prediction error in Table \ref{tab:shocksL1}. 	
The individual and averaged errors are larger for this example, but again the errors we get are in line with the expectations from the original wPINN experiments.

The rest of the discussion is analogous to the standing shock. However, because the shock is moving at a constant velocity, the PDE residual now essentially amounts to learning a single linear transformation to ensure correct placement of the solution profile learned from the initial data. This is slightly less straightforward than previously, resulting in a larger relative error and a small advantage for the modified wPINNs.

\subsection{Rarefaction Wave}
As another prototypical example, we consider a rarefaction wave on $[0,0.5)\times[-1,1]$ for the initial data 
\begin{equation}
\soln_0(x) = \begin{cases}	-1 \quad&\text{for } x \leq 0, \\ 1 \quad& \text{for } x > 0.	\end{cases}
\end{equation}
For this, the exact solution is given by a rarefaction wave 
\begin{equation}
\soln(t,x) = \begin{cases}	-1 \quad&\text{for } x \leq -t, \\ \frac{x}{t} \quad&\text{for } -t < x \leq t, \\ 1 \quad& \text{for } x > t.	\end{cases}
\end{equation}
Without additional entropy conditions, weak solutions are non-unique for the given initial data. However, due to the entropy residual, the modified wPINNs find the correct entropy solution. We train the neural networks for this example using the same hyperparameters as in the previous examples, except that we train for a total of $N_\mathrm{ep} = 3000$ epochs now.
\begin{figure}[ht]
	\begin{subfigure}[t]{0.45\textwidth}
		\vskip 0pt
		\includegraphics[width=\textwidth]{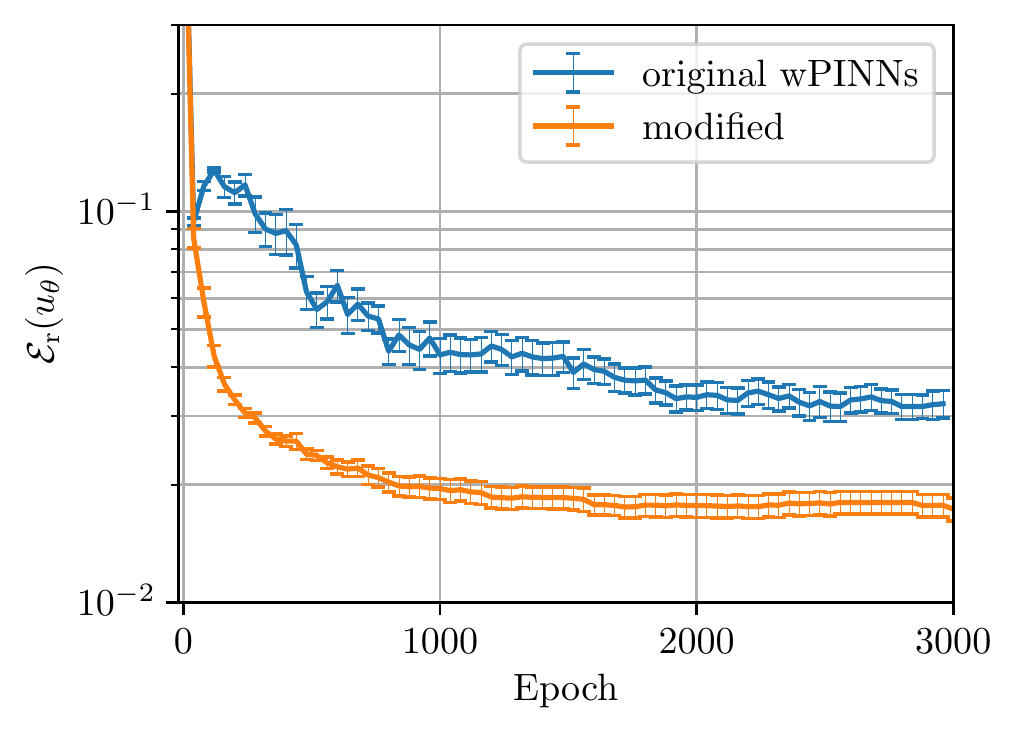}
		\caption{Evolution of the mean value and standard error of $\relerr$ for the respective network ensembles. After 800 epochs, the modified wPINNs are already within $15\%$ of the final relative error, while the original wPINNs are still $35\%$ worse than their final relative error.}
		\label{fig:rarefacL1}
	\end{subfigure}
	\hspace{1em}
	\begin{subfigure}[t]{0.45\textwidth}
		\vskip 0pt
		\includegraphics[width=\textwidth]{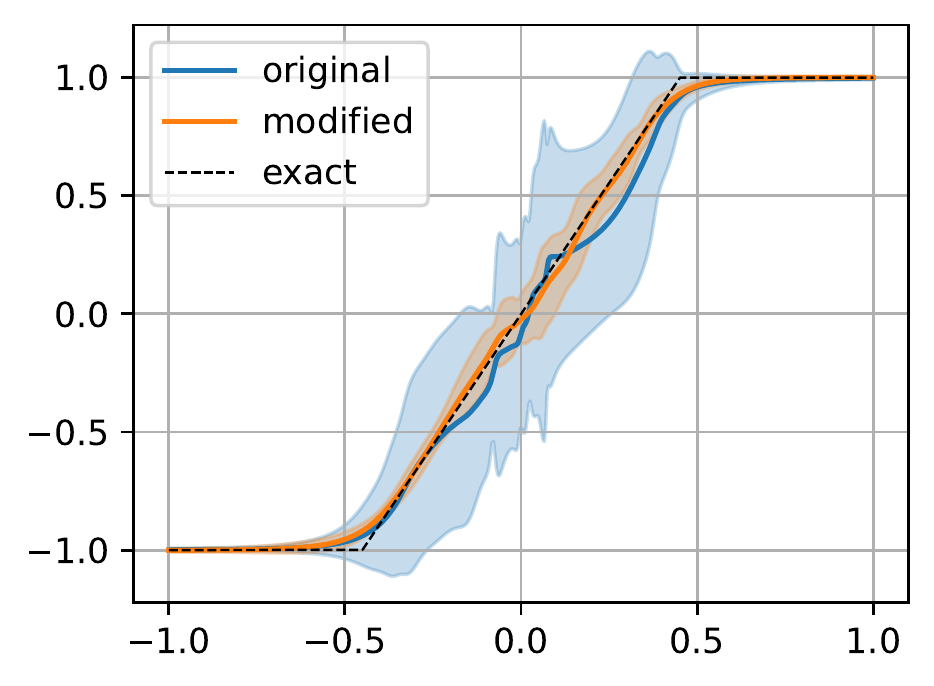}
		\caption{Comparison of average network predictions for original and modified wPINNs at  $t = 0.45$ after $ N_{\mathrm{ep}} = 800$ epochs. Shaded bands are $2\sigma$-standard deviations of the average network prediction.	We refer to the text for discussion.}
		\label{fig:rarefacSlices}
	\end{subfigure}
	\caption{Rarefaction wave with shock initial data}
\end{figure}

The mean  relative error of the ensemble as training progresses is shown in Figure \ref{fig:rarefacL1}. Because this example is no longer dominated by how quickly the network can learn the correct shape of the solution from the initial data, it is better-suited to illustrate the advantages of the modified wPINNs. In this example we see that the modified wPINNs provide a sizable advantage over the original wPINN algorithm.

For the errors of the original and modified wPINNs average network predictions, we refer to Table \ref{tab:shocksL1}. { After 800 epochs of training the modified wPINNs are already within $15\%$ of their final accuracy. At this epoch, the modified wPINNs are roughly twice as accurate as the original wPINNs. The original wPINN network predictions have, on average, a {$4.3\% \pm 0.42\%$} relative error while the modified wPINNs have an error of only {$2.0\% \pm 0.11\%$}. Prolonged training is able to close the accuracy gap somewhat, as the modified wPINNs are closer to finishing learning and additional epochs enable the original wPINNs to catch up slightly. However, they are not able to completely close the accuracy gap to the modified wPINNs even after 3.5x longer training.}

In Figure \ref{fig:rarefacSlices} we show the ensemble averaged prediction   after 800 epochs of training for the original and modified wPINN algorithm at $t = 0.45$.  We see that the accuracy of the original wPINNs is much lower, such that one can easily see the deviations from the reference solution.
Further, the original wPINN algorithm produces spurious oscillations when approximating the solution to this problem. While both the modified loss and the original wPINN loss penalize spurious oscillations as argued in Section \ref{sec:hypCons}, we observe that the modified approach seems to penalize this more effectively, so we do not observe significant oscillations for its results. 

Finally, we observe that the standard deviation of the average network prediction is much smaller for the modified wPINNs, showing that the modified wPINNs are more stable during training. This is an advantage, because  every network in our ensemble is an accurate approximation. The \emph{worst-performing} network in the modified wPINNs ensemble at epoch 800 achieves a relative $L^1$ error of $2.91\%$, which is as good as the \emph{averaged} network prediction of the original wPINNs at $2.94\%$ as shown in Table \ref{tab:shocksL1}. Thus, the smaller standard deviations for the modified wPINNs mean that one could also forego ensemble training (or use smaller ensembles), reducing computational costs.

\subsection{Sine Wave}
Lastly we consider the Burgers equation for sine initial data $u_0(x) = -\sin(\pi x)$ on $[0,1)\times[-1,1]$ and $u(t,-1) = u(t,1) = 0$ for $t \in [0,1)$. This example is significantly more challenging to solve to high precision using neural networks, because unlike previous examples it contains shock formation from smooth initial data.  This means the network approximation has to steepen over time, developing into a shock separating two rarefactions. This example is also studied in \cite{wPINNS} and represents the typical performance of both original and modified wPINNs better than the more synthetic tests of a single discontinuity or a pure rarefaction example.

Because this example is more challenging, and to reproduce the simulation setup in \cite{wPINNS}, instead of uniformly randomly sampled points on the domain we use low-discrepancy Sobol points. These approximate the numerical integrals more precisely than previous uniformly randomly sampled points, improving network performance.

In our preliminary experiments we find the hyperparameter choices from the previous two examples to perform best for the original wPINN algorithm again. However, the exact solution is more complicated, so we increase the number of epochs to $N_\mathrm{ep} = 75000$ to match the training outlined in \cite{wPINNS} for this example. However, we find that modified wPINNs require significantly less epochs to achieve good results during our numerical experiments.  Because we train for more epochs, we reduce $\gamma$ to $0.015$ for the exponential loss averages.

\begin{figure}[ht]
	\centering
	\begin{subfigure}[t]{0.45\textwidth}
		\vskip 0pt
		\includegraphics[width=\textwidth]{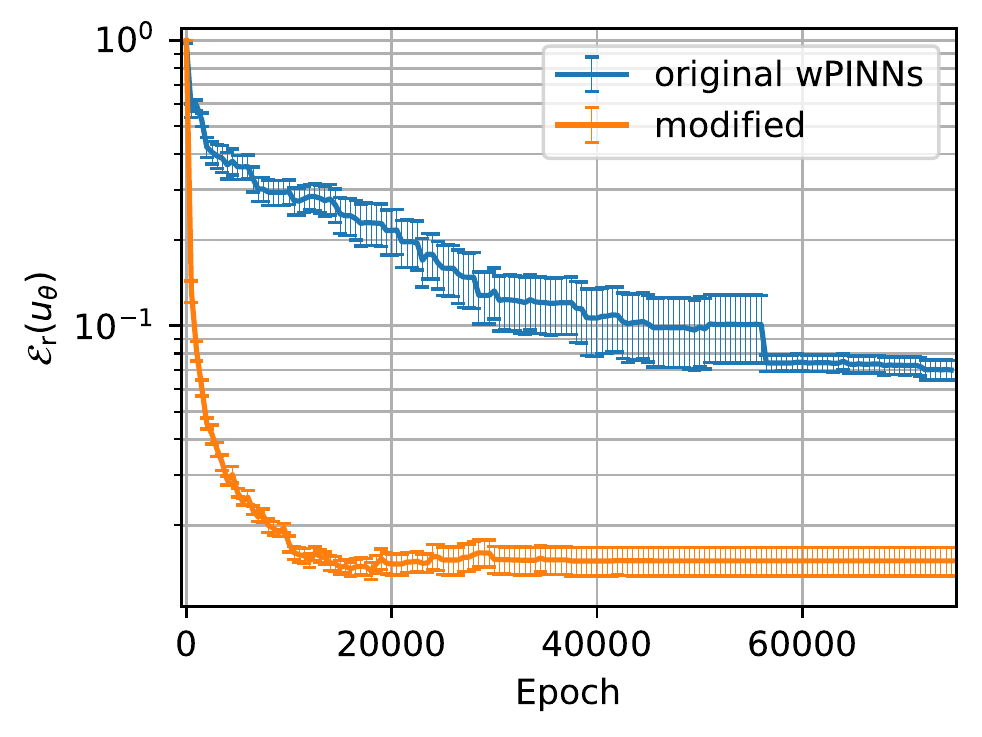}
		\caption{Evolution of the mean value and standard error of $\relerr$ for the respective network ensembles. The modified wPINNs are within $15\%$ of their final accuracy after 10000 epochs, while the original wPINNs are still four times less accurate than they are at the end of training.}
		\label{fig:sineL1}
	\end{subfigure}	\hspace{1em}
	\begin{subfigure}[t]{0.45\textwidth}
		\vskip 0pt
		\includegraphics[width=\textwidth]{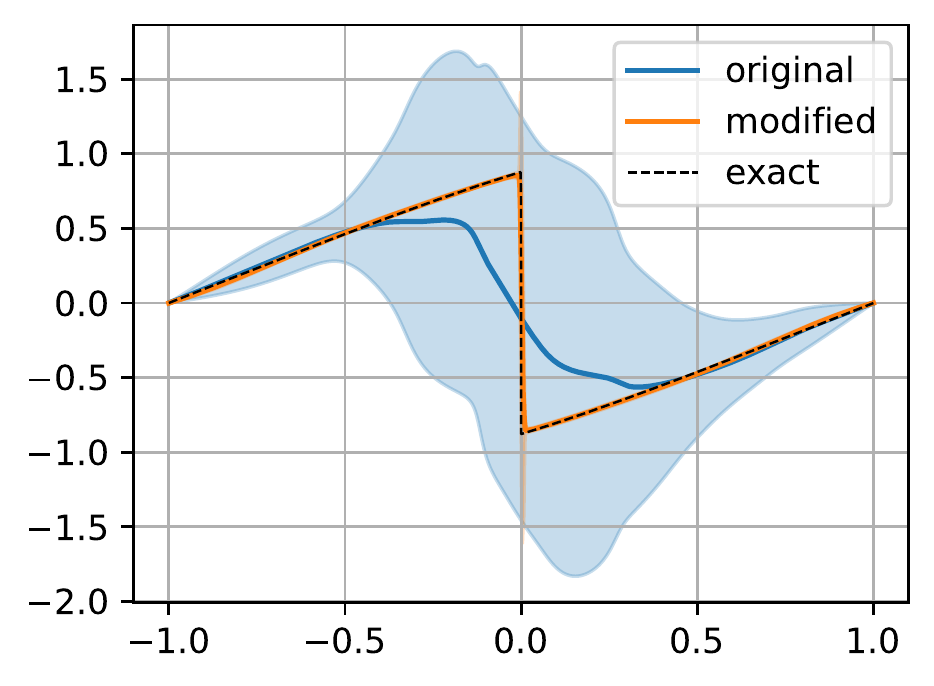}
		
		\caption{Comparison of average network predictions for original and modified wPINNs at  $t = 0.75$ after $N_{\mathrm{ep}} = 10000$ epochs. Shaded bands are $2\sigma$-standard deviations of the average network prediction. We refer to the text for discussion.}
		\label{fig:sineSlices}
	\end{subfigure}
	\caption{Solution and training progress for sine initial data}
\end{figure}

As before, we show the mean relative error of the ensemble over the training epochs in Figure \ref{fig:sineL1}. In this example we see a very noticeable difference between the original and the modified wPINN algorithms. The mean value of $\relerr$ at the final training epoch for the original wPINNs  is $7.0\% \pm 0.6\%$, while the modified algorithm gives a mean relative error of only $1.5\% \pm 0.17 \%$. We also compare results after 10000 epochs, when the modified wPINNs are within $15\%$ of their final accuracy. The modified wPINNs achieve an accuracy of $1.73\% \pm 0.11 \%$, while the original wPINNs have an accuracy of only $29.5\% \pm 2.9\%$ on average.  The drop in the mean relative error for the original wPINNs after about $55000$ epochs alongside a large reduction in the standard error is due to a single poorly-performing outlier network improving to match the performance of the remaining ensemble at this point in training. This is very visible given the limited ensemble size of $N^* = 16$ retrainings.

The average network predictions after 10000 epochs are compared in Figure \ref{fig:sineSlices}. { At this time during training, the original wPINNs have a  relative error of $\relerr(\soln^*_\mathrm{orig}) = 21.4 \%$, while the modified wPINNs already achieve an accuracy of $\relerr(\soln^*_\mathrm{mod}) = 1.53 \%$ }. At the end of training, the original wPINNs have a  relative error of $\relerr(\soln^*_\mathrm{orig}) = 5.03 \%$, which is larger but comparable to the results from \cite{wPINNS}. For the modified wPINNs we find $\relerr(\soln^*_\mathrm{mod}) = 1.19\ \%$ instead, which is better than the results from \cite{wPINNS}. { The modified wPINNs have finished learning completely after about 20000 epochs, as their relative error plateaus. Any prolonged training benefits only the original wPINNs. Yet, even after the full 75000 epochs, the original wPINNs are not as accurate as the modified wPINNs.} 
Again, the modified wPINN ensemble has a much smaller standard deviation than the original wPINN ensemble. This behavior persists even after the full training period of 75000 epochs and is even more pronounced than in the previous example.

Note that for the modified wPINNs algorithm, its final level of accuracy is roughly in line with the previous example, with the main difference being that it takes more epochs to reach this level of accuracy. In our tests we saw that increasing the number of collocation points and choosing a different learning rate schedule can further improve results, while the networks themselves appear to be able to approximate the solution of the underlying problem to higher accuracy without changing the network size, indicating that these hyperparameters restrict the final precision most here.
	\section{Numerical Results - Extensions}
\label{sec:extensions}
\subsection{Weak Boundary Conditions}
\label{sec:weakBC_num}
We test our extension for weak Dirichlet boundary data. We solve the Burgers equation on $[0,0.8) \times [-0.5,0.5]$ for the initial data
\begin{equation}
u_0 = \begin{cases}	2 \quad&\text{for } x \leq 0, \\ 0 \quad& \text{for } x > 0,	\end{cases}
\label{eqn:weakBC_icond}
\end{equation}
and the boundary data $g(t,-0.5) = 2 - 2t$ and $g(t,0.5) = -1$ for all $t\in [0,0.8)$. In this case, the left-side boundary is part of the inflow boundary the entire time, while the right-side boundary initially starts as inflow boundary and then switches to outflow boundary as the shock propagates out of the domain.

For the training of this problem, as described in section \ref{sec:weakBC_theo}, we subdivide the time domain into six equally large time slices because long time domains can be problematic when using weak boundary conditions. Then we begin by training a neural network on the first time slice $t \in [0,2/15]$ and use the evaluation of the first neural network at $t = {2}/{15}$ as the initial data for a separate neural network trained on $t \in [{2}/{15},{4}/{15}]$ and repeat this process until the final time is reached.

For training we use the following hyperparameters: $\tau_{\mathrm{min}} = 0.005$, $\tau_{\mathrm{max}} = 0.009$, $N_\mathrm{max} = 8$, $N_\mathrm{ep} = 5000$. The network $\soln_\theta$ has $l = 6$ layers with a width of $w = 20$ and $\tanh$-activation, while the networks $\varphi_\chi$ and $\xi_\nu$ use $4$ layers of width $14$ and $\tanh$-activation. We use $N_\mathrm{int} = 32768$ points per time slice and $N_\mathrm{ic} = N_\mathrm{bc} = 4096$ points for initial- and boundary conditions on each subdomain.  We do not perform any ensemble training. The hyperparameters were not exhaustively tuned to achieve best performance.
\begin{figure}[h]
	\centering
	\begin{subfigure}[t]{0.45\textwidth}
		\vskip 0pt
		\includegraphics[width=\textwidth]{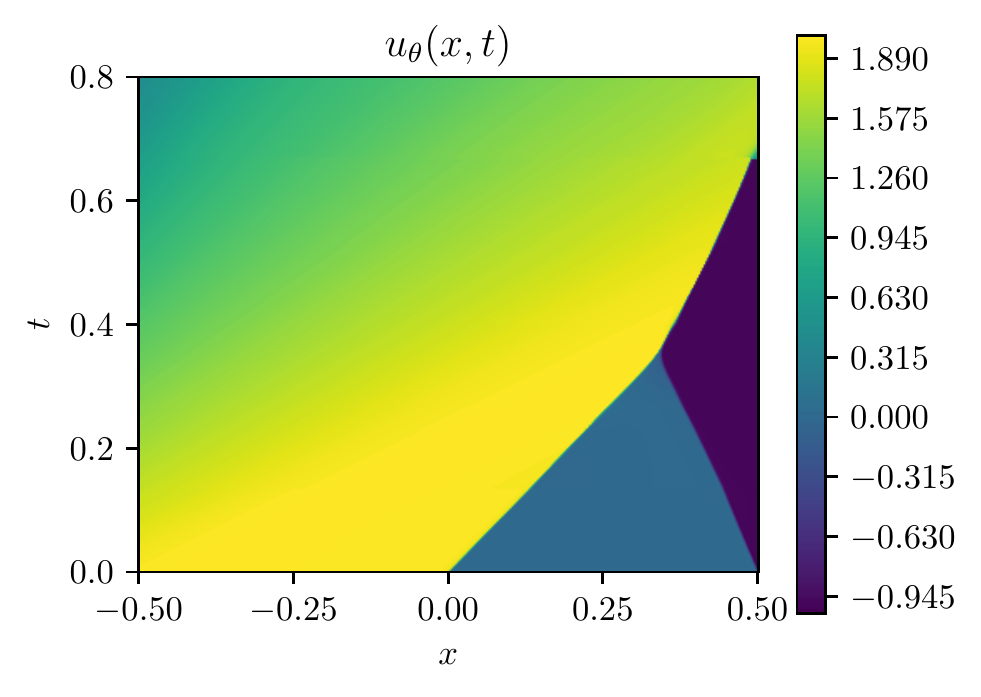}
		\caption{Neural network solution. The network automatically decide when the right boundary is part of the inflow- or outflow boundary.}
		\label{fig:weakBC}
	\end{subfigure}	\hspace{1em}
	\begin{subfigure}[t]{0.49\textwidth}
		\vskip 0pt
		\includegraphics[width=\textwidth]{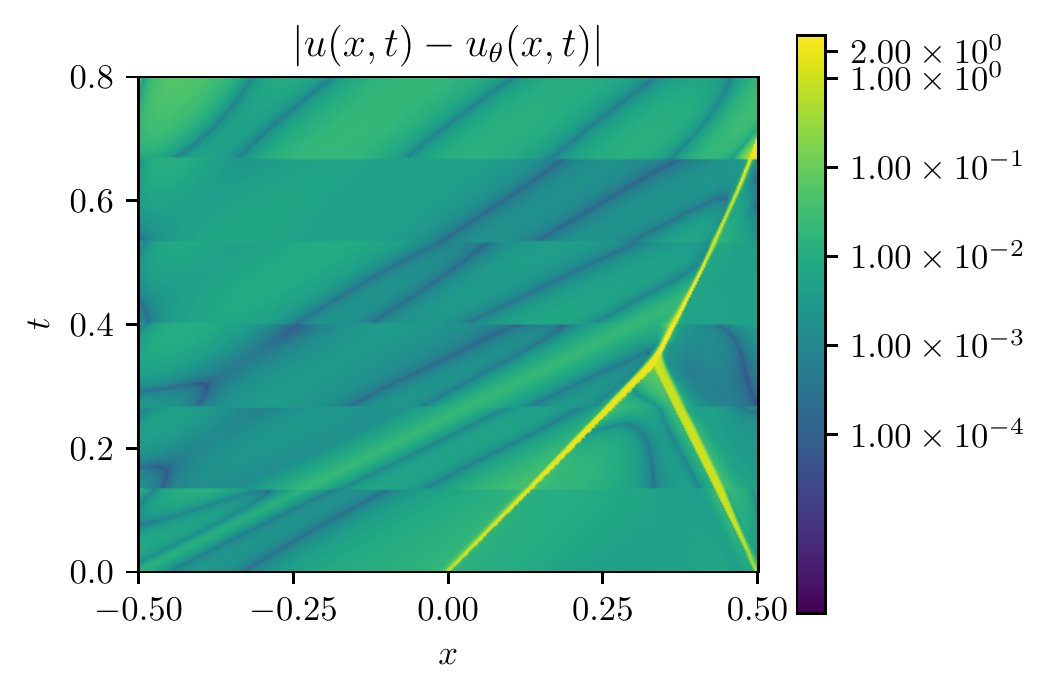}			
		\caption{Pointwise error of the neural network approximation. The logarithmic scale of the colorbar makes the small discontinuities in the pointwise error at the interfaces between the time slices visible.}
		\label{fig:weakBC_diff}
	\end{subfigure}
	\caption{Burgers equation with weak boundary conditions.}
\end{figure}

The neural network approximation after training is displayed in Figure \ref{fig:weakBC} and the pointwise error to the exact solution in Figure \ref{fig:weakBC_diff}.

We see that both imposing the boundary condition when valid, as well as \emph{not} enforcing the boundary data on the right as the shock exits the domain, work correctly. In particular, the shock exits the domain around $t = 0.7$. Note that due to understanding the boundary conditions in a weak sense this does not require any prior knowledge but rather results from simulation. The logarithmically scaled pointwise error shows the discontinuities at the interface between two time slices clearly and that coupling the subdomains does not incur large errors. Our approach solves the problem with good accuracy, leading to a final relative $L^1$-error of $\relerr(u_\theta) = 1.21\%$ across the entire domain.
\subsection{Compressible Euler Equations}
\label{sec:numResults_Sod}
In this section we apply the modified wPINNs to the compressible Euler equations for the Sod shock tube.
The compressible Euler equations in primitive variables are given by 
\begin{empheq}[left=\empheqlbrace]{align}
\rho_t + (\rho v)_x &= 0, \label{eqn:euler1}\\
(\rho v)_t + (\rho v^2 + p)_x &= 0\label{eqn:euler2},\\
E_t + \left(v (E + p)\right)_x &= 0\label{eqn:euler3},
\end{empheq} $\text{in } \tDom \times \symbb{R}$
where $\rho$ is the density, $v$ is the velocity and $p$ is the pressure. For Sod's problem, the pressure relates to the internal energy $e$ of the gas as $p = \rho e (\gamma - 1)$ with $\gamma$ the adiabatic exponent of the gas. We consider $\gamma=1.4$, as it would be for an ideal diatomic gas. Lastly, the total energy $E$ is given by $E = \rho e + \frac{1}{2} \rho v^2$. As in \cite{STRELOW2023112041}, our network parameterizes the primitive variables $(\rho,p,v)$ instead of conservative variables $(\rho,\rho v, E)$ because this avoids divisions by zero when computing $v = (\rho v)/\rho$. However, the system is still given in conservative form, to align with our previous arguments in favor of using weak norms.

Sod's problem then consists of solving the Euler equations for the Riemann problem with initial data 
\begin{empheq}[left=\empheqlbrace]{align}
(\rho,p,v) &= (1.0,1.0,0.0)  &&\text{for } t=0,\,\; x < 0.5 \text{ and}\\
(\rho,p,v) &= (0.125,0.1,0.0) &&\text{for }t=0,\,\;  x \geq 0.5.
\end{empheq}

An exact solution is available for comparison with the numerical results \cite{SOD19781}.

We parameterize our approximations for $\rho$, $v$ and $p$ each using a separate neural network. While one could also use a single neural network to describe the entire state vector $(\rho,v,p)$, this holds the disadvantage that discontinuities in one state variable tend to lead to poor performance in the other state variables, especially for smaller networks, because all components share the same overall network. This makes learning significantly more difficult. Likewise, we have three separate networks responsible for the weak norm estimation, one per equation of the Euler system.

For training, we use a PDE residual loss akin to equation \eqref{eqn:mixedRes} employing integration by parts on the $x$-derivatives for each of \eqref{eqn:euler1} -- \eqref{eqn:euler3}, and sum up the contributions per equation. We use $q = 2$ in the modified loss and this seems to work well.

Lastly, for the entropy residual part in our approach we use the entropy-entropy flux pair 
\begin{equation}
(\eta,q) \coloneqq (\rho S, \rho v S),
\end{equation}
with $S \coloneqq \ln(\rho^\gamma/p)$ and enforce the usual entropy inequality \cite{Svaerd2016}. Initial- and boundary conditions are incorporated into the loss through their standard $L^2$ contributions, completely analogously to the scalar case. We do not perform any ensemble training for this example.
\begin{figure}[t]
	\centering\includegraphics[width = 0.6\textwidth]{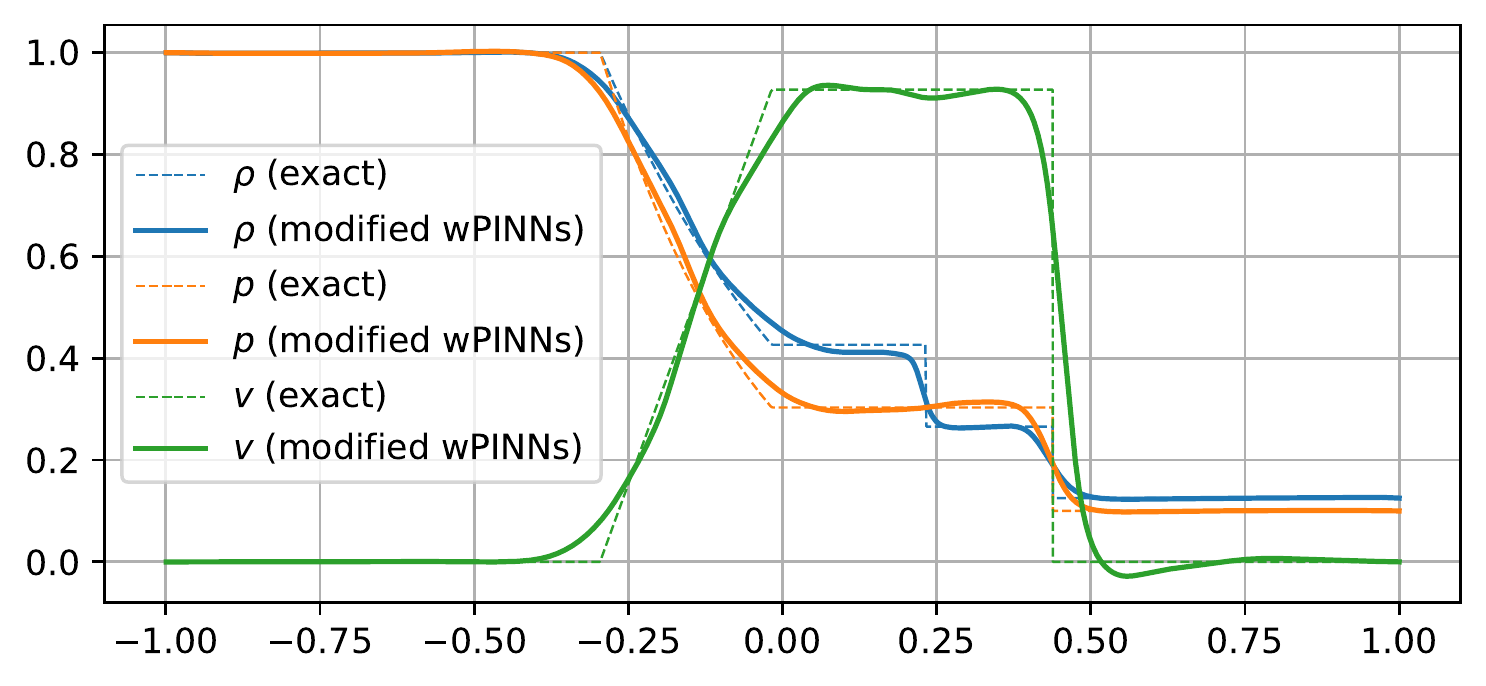}
	\caption{Density, pressure and velocity of the neural network approximation for Sods shock tube at $t = 0.25$. We see that all the different regions of the exact solution are resolved and there is only mild smearing around the different discontinuities.}
	\label{fig:sodSlice}
\end{figure}	
Our neural networks for the approximation of $\rho,v$ and $p$ consist of $6$ fully connected layers of 45 neurons each. This is larger than in our previous scalar examples, but we find that this size is necessary to adequately capture the two discontinuities that will arise in the density profile. As activation function we use a scaled $\operatorname{softplus}$ function given by $\operatorname{softplus}(40; x) \coloneqq \frac{1}{40}\log(1 + \exp(40\cdot x))$, which serves as a smooth approximation to the $\operatorname{ReLu}$-function. The scaling factor $40$ was chosen arbitrarily. We find that this smooth approximation performs better than the $\operatorname{ReLu}$ activation function, presumably because having continuous derivatives makes learning easier. We apply the unscaled $\operatorname{softplus}$-function to the output of our density- and pressure networks. The entropy-entropy flux pair is not well-defined for negative density or pressure, however basic random initialization of the network weights cannot ensure that this does not occur. This step then guarantees that the output of these networks is always positive and is important for the early steps during training. Later into the training, the approximation is reasonably close to the true solution, so it has non-negative pressure and density anyway.

All other networks in this example have 6 layers with 24 neurons each and also use ${\operatorname{softplus}(40 ; x)}$ as layer activation function.

We solve the shock tube on $[0,0.4]\times [-1,1]$. We use $N_{\mathrm{int}} = 65536$ Sobol points in the interior of the domain, and $N_{\mathrm{ic}} = N_{\mathrm{bc}} = 8096$ points for initial and boundary condition respectively. We train for a total of $N_{\mathrm{ep}} = 42000$ epochs using a learning rate of $\tau_{\mathrm{min}} = 0.005$ and $\tau_{\mathrm{max}} = 0.014$. We choose a learning rate schedule that linearly reduces the learning rate by a factor of $30$ over the course of training, because we need to minimize the loss to fairly high precision. We optimize using the Adam optimizer, and set the hyperparameters $\beta_1 = 0.2$ and $\beta_2 = 0.99$. Using lower values of $\beta_1$ and $\beta_2$ than “default” is common for training adversarial networks \cite{gemp2019}.

We show the final trained network at time $t = 0.25$ in Figure \ref{fig:sodSlice}, which is representative of the performance over the entire space-time domain. We compute the relative $L^1$ in space-time error for the density, pressure and velocity respectively and find $\relerr(\rho_\mathrm{mod}) = 1.38 \%$, $\relerr(p_\mathrm{mod}) = 1.47 \%$ and $\relerr(v_\mathrm{mod}) = 3.77 \%$. 

	\section{Discussion}
	 Approximating nonlinear hyperbolic conservation laws  using conventional PINNs fails when the exact solution is discontinuous and PDE residuals are minimized in the $L^2$ norm. The failure of conventional PINNs is structural and cannot be remedied by using larger networks or training more, because the residual is measured in an unsuitable norm. This class of problems requires using \emph{weak norms} instead. The wPINN strategy from \cite{wPINNS} enables the accurate approximation of entropy solutions to scalar conservation laws by computing weak norms. Test functions for the PDE residual are parameterized with neural networks, leading to an adversarial min-max problem structure which may be solved with standard deep learning techniques.
	 
	 However, the training of wPINNs is computationally more expensive than that of conventional PINNs, because of the min-max problem structure. Using the supremum-based definition of the weak norms leads to a challenging maximization problem.
	 
	 We improve on this by providing a new strategy for approximating weak norms. The maximizer  can be found as the solution to a dual elliptic problem. This avoids normalizing the test function. While this still leads to a min-max formulation of the loss, it has a simpler structure that makes training wPINNs faster.
	 
	 Beyond the changes to computing dual norms, we modify several core parts of the original wPINN algorithm to accelerate training and enable extensions to systems of hyperbolic conservation laws. In the original wPINN algorithm, the Kruzhkov entropy residual has two important roles: enforcing the PDE in a weak sense and selecting the entropy solution. We split these roles into separate terms by introducing a PDE loss and an entropy loss with respect to a strictly convex entropy-entropy flux pair. This gives more flexibility to both terms. Crucially, this means we do not have to use the Kruzhkov family of entropies and can opt for a smoother entropy, which makes training easier. An additional benefit to not using the Kruzhkov entropies is that our modified wPINNs naturally extend to systems of hyperbolic conservation laws endowed with a convex entropy.
	 
	 
	 Numerical experiments illustrate the benefits of the modified wPINNs. Especially for more challenging problems we see large improvements expressed by lower mean relative errors after the same number of training epochs and better approximations at the end of training. The modified wPINNs finish training in significantly less epochs than the original wPINNs. 
	 
	 
	 Finally, we show that the modified wPINNs are suitable for weak boundary conditions and systems of hyperbolic conservation laws. For scalar conservation laws we cover weak boundary conditions by enforcing an entropy-based inequality on the boundary of the domain. We verify that this approach is effective by solving the Burgers equation with weak boundary conditions. Our experiments show that this approach enables shocks to exit the computational domain without requiring any prior knowledge of when this occurs.	 
	 Further, we give an example of modified wPINNs naturally being applicable to systems of hyperbolic conservation laws, by summing over the PDE residuals of each separate equation. We solve Sod's shock tube problem for the Euler equations successfully.

\section*{Acknowledgements}	
The authors gratefully acknowledge the computing time provided to them on the high-performance computer Lichtenberg at the NHR Centers NHR4CES at TU Darmstadt. This is funded by the Federal Ministry of Education and Research, and the state governments participating on the basis of the resolutions of the GWK for national high performance computing at universities.

\printbibliography

@article{wPINNS,
  doi = {10.48550/ARXIV.2207.08483},
  
  url = {https://arxiv.org/abs/2207.08483},
  
  author = {De Ryck, Tim and Mishra, Siddhartha and Molinaro, Roberto},
  
  title = {wPINNs: Weak Physics informed neural networks for approximating entropy solutions of hyperbolic conservation laws},
  
  publisher = {arXiv},
  
  year = {2022},
  
  copyright = {arXiv.org perpetual, non-exclusive license},
  archivePrefix={arXiv},
  primaryClass={math.NA},
  journal={arXiv:2207.08483}
}

@article{Lye_2020,
	doi = {10.1016/j.jcp.2020.109339},
  
	url = {https://doi.org/10.1016%2Fj.jcp.2020.109339},
  
	year = 2020,
	month = {jun},
  
	publisher = {Elsevier {BV}
},
  
	volume = {410},
  
	pages = {109339},
  
	author = {Kjetil O. Lye and Siddhartha Mishra and Deep Ray},
  
	title = {Deep learning observables in computational fluid dynamics},
  
	journal = {Journal of Computational Physics}
}

@article{panov1994uniqueness,
  title={Uniqueness of the solution of the Cauchy problem for a first order quasilinear equation with one admissible strictly convex entropy},
  author={Panov, E. Yu.},
  journal={Mathematical Notes},
  volume={55},
  number={5},
  pages={517--525},
  year={1994},
  publisher={Kluwer Academic Publishers-Plenum Publishers}
}

@article{de2004minimal,
  title={Minimal entropy conditions for Burgers equation},
  author={De Lellis, Camillo and Otto, Felix and Westdickenberg, Michael},
  journal={Quarterly of applied mathematics},
  volume={62},
  number={4},
  pages={687--700},
  year={2004}
}

@article{Chen_2020,
	doi = {10.4208/cmr.2020-0051},
  
	url = {https://doi.org/10.4208%2Fcmr.2020-0051},
  
	year = 2020,
	month = {jun},
  
	publisher = {Global Science Press},
  
	volume = {36},
  
	number = {3},
  
	pages = {354--376},
  
	author = {Jingrun Chen},
  
	title = {A Comparison Study of Deep Galerkin Method and Deep Ritz Method for Elliptic Problems with Different Boundary Conditions},
  
	journal = {Communications in Mathematical Research}
}

@inproceedings{kingma_Adam,
  author    = {Diederik P. Kingma and
               Jimmy Ba},
  editor    = {Yoshua Bengio and
               Yann LeCun},
  title     = {Adam: {A} Method for Stochastic Optimization},
  booktitle = {3rd International Conference on Learning Representations, {ICLR} 2015,
               San Diego, CA, USA, May 7-9, 2015, Conference Track Proceedings},
  year      = {2015},
  url       = {http://arxiv.org/abs/1412.6980},
  bibsource = {dblp computer science bibliography, https://dblp.org}
}

@inproceedings{reddi_amsgrad,
  author    = {Sashank J. Reddi and
               Satyen Kale and
               Sanjiv Kumar},
  title     = {On the Convergence of Adam and Beyond},
  booktitle = {6th International Conference on Learning Representations, {ICLR} 2018,
               Vancouver, BC, Canada, April 30 - May 3, 2018, Conference Track Proceedings},
  IGNOREpublisher = {OpenReview.net},
  year      = {2018},
  url       = {https://openreview.net/forum?id=ryQu7f-RZ},
  bibsource = {dblp computer science bibliography, https://dblp.org}
}

@INPROCEEDINGS{smith_cyclical,

  author={Smith, Leslie N.},

  booktitle={2017 IEEE Winter Conference on Applications of Computer Vision (WACV)}, 

  title={Cyclical Learning Rates for Training Neural Networks}, 

  year={2017},

  volume={},

  number={},

  pages={464-472},

  doi={10.1109/WACV.2017.58}}

@Article{E2018,
author={E, Weinan
and Yu, Bing},
title={The Deep Ritz Method: A Deep Learning-Based Numerical Algorithm for Solving Variational Problems},
journal={Communications in Mathematics and Statistics},
year={2018},
month={Mar},
day={01},
volume={6},
number={1},
pages={1-12},
issn={2194-671X},
doi={10.1007/s40304-018-0127-z},
url={https://doi.org/10.1007/s40304-018-0127-z}
}

@article{pyclaw-sisc,
        Author = {Ketcheson, David I. and Mandli, Kyle T. and Ahmadia, Aron J. and Alghamdi, Amal and {Quezada de Luna}, Manuel and Parsani, Matteo and Knepley, Matthew G. and Emmett, Matthew},
        Journal = {SIAM Journal on Scientific Computing},
        Month = nov,
        Number = {4},
        Pages = {C210--C231},
        Title = {{PyClaw: Accessible, Extensible, Scalable Tools for Wave Propagation Problems}},
        Volume = {34},
        Year = {2012}}

@article{BLN_boundary,
author = {   C.   Bardos  and    A. Y.   Leroux  and    J. C.   Nedelec },
title = {First order quasilinear equations with boundary conditions},
journal = {Communications in Partial Differential Equations},
volume = {4},
number = {9},
pages = {1017-1034},
year  = {1979},
publisher = {Taylor & Francis},
doi = {10.1080/03605307908820117},

URL = { 
        https://doi.org/10.1080/03605307908820117
    
},
eprint = { 
        https://doi.org/10.1080/03605307908820117
    
}
}

@article {kondo2001measure,
    AUTHOR = {Kondo, C. I. and LeFloch, P. G.},
     TITLE = {Measure-valued solutions and well-posedness of
              multi-dimensional conservation laws in a bounded domain},
   JOURNAL = {Port. Math. (N.S.)},
  FJOURNAL = {Portugaliae Mathematica. Nova S\'{e}rie},
    VOLUME = {58},
      YEAR = {2001},
    NUMBER = {2},
     PAGES = {171--193},
      ISSN = {0032-5155},
   MRCLASS = {35L65 (35B30 76N10)},
  MRNUMBER = {1836261},
MRREVIEWER = {Hermano Frid},
}

@book {dafermos,
    AUTHOR = {Dafermos, Constantine M.},
     TITLE = {Hyperbolic conservation laws in continuum physics},
    SERIES = {Grundlehren der mathematischen Wissenschaften [Fundamental
              Principles of Mathematical Sciences]},
    VOLUME = {325},
   EDITION = {Fourth},
 PUBLISHER = {Springer-Verlag, Berlin},
      YEAR = {2016},
     PAGES_IGNORE = {xxxviii+826},
      ISBN_IGNORE = {978-3-662-49449-3; 978-3-662-49451-6},
   MRCLASS = {35L65 (35-02 35L67 74J40 76-02)},
  MRNUMBER = {3468916},
MRREVIEWER = {Marta Lewicka},
       DOI = {10.1007/978-3-662-49451-6},
       URL = {https://doi.org/10.1007/978-3-662-49451-6},
}

@article{SirignanoDGM,
title = {DGM: A deep learning algorithm for solving partial differential equations},
journal = {Journal of Computational Physics},
volume = {375},
pages = {1339-1364},
year = {2018},
issn_ignore = {0021-9991},
doi = {https://doi.org/10.1016/j.jcp.2018.08.029},
url = {https://www.sciencedirect.com/science/article/pii/S0021999118305527},
author = {Justin Sirignano and Konstantinos Spiliopoulos}
}

@article{raissi2019physics,
  title={Physics-informed neural networks: A deep learning framework for solving forward and inverse problems involving nonlinear partial differential equations},
  author={Raissi, Maziar and Perdikaris, Paris and Karniadakis, George E},
  journal={Journal of Computational Physics},
  volume={378},
  pages={686--707},
  year={2019},
  publisher={Elsevier}
}

@article{liu2022discontinuity,
  title={Discontinuity Computing with Physics-Informed Neural Network},
  author={Liu, Li and Liu, Shengping and Yong, Heng and Xiong, Fansheng and Yu, Tengchao},
  journal={arXiv:2206.03864},
  year={2022}
}

@article{DERYCK2021732,
title = {On the approximation of functions by tanh neural networks},
journal = {Neural Networks},
volume = {143},
pages = {732-750},
year = {2021},
issn = {0893-6080},
doi = {https://doi.org/10.1016/j.neunet.2021.08.015},
url = {https://www.sciencedirect.com/science/article/pii/S0893608021003208},
author = {Tim {De Ryck} and Samuel Lanthaler and Siddhartha Mishra}
}

@article{mishra_L2fail,
    author = {Mishra, Siddhartha and Molinaro, Roberto},
    title = "{Estimates on the generalization error of physics-informed neural networks for approximating PDEs}",
    journal = {IMA Journal of Numerical Analysis},
    year = {2022},
    month = {01},
    issn = {0272-4979},
    doi = {10.1093/imanum/drab093},
    url = {https://doi.org/10.1093/imanum/drab093},
    note = {drab093},
    eprint = {https://academic.oup.com/imajna/advance-article-pdf/doi/10.1093/imanum/drab093/42198350/drab093.pdf},
}

@article{PinnOrig1,
author = {Dissanayake, M. and Phan-Thien, N.},
title = {Neural-network-based approximations for solving partial differential equations},
journal = {Communications in Numerical Methods in Engineering},
volume = {10},
number = {3},
pages = {195-201},
doi = {https://doi.org/10.1002/cnm.1640100303},
url = {https://onlinelibrary.wiley.com/doi/abs/10.1002/cnm.1640100303},
eprint = {https://onlinelibrary.wiley.com/doi/pdf/10.1002/cnm.1640100303},
year = {1994}
}

@ARTICLE{PinnOrig2,

  author={Lagaris, I.E. and Likas, A.C. and Papageorgiou, D.G.},

  journal={IEEE Transactions on Neural Networks}, 

  title={Neural-network methods for boundary value problems with irregular boundaries}, 

  year={2000},

  volume={11},

  number={5},

  pages={1041-1049},

  doi={10.1109/72.870037}}

@ARTICLE{PinnOrig3,

  author={Lagaris, I.E. and Likas, A. and Fotiadis, D.I.},

  journal={IEEE Transactions on Neural Networks}, 

  title={Artificial neural networks for solving ordinary and partial differential equations}, 

  year={1998},

  volume={9},

  number={5},

  pages={987-1000},

  doi={10.1109/72.712178}}

@article{variPinn,
  doi = {10.48550/ARXIV.1912.00873},
  
  url = {https://arxiv.org/abs/1912.00873},
  
  author = {Kharazmi, E. and Zhang, Z. and Karniadakis, G. E.},
  
  title = {Variational Physics-Informed Neural Networks For Solving Partial Differential Equations},
  
  publisher = {arXiv},
  
  year = {2019},
  journal={arXiv:1912.00873},
  copyright = {arXiv.org perpetual, non-exclusive license}
}

@article {Guermond_L1,
    AUTHOR = {Guermond, Jean-Luc and Marpeau, Fabien and Popov, Bojan},
     TITLE = {A fast algorithm for solving first-order {PDE}s by
              {$L^1$}-minimization},
   JOURNAL = {Commun. Math. Sci.},
  FJOURNAL = {Communications in Mathematical Sciences},
    VOLUME = {6},
      YEAR = {2008},
    NUMBER = {1},
     PAGES = {199--216},
      ISSN = {1539-6746},
   MRCLASS = {65N35 (35J05 35R25)},
  MRNUMBER = {2398004},
MRREVIEWER = {Monika Neda},
       URL = {http://projecteuclid.org/euclid.cms/1204905784},
}

@article{fractionalADE,
author = {Pang, Guofei and Lu, Lu and Karniadakis, George Em},
title = {fPINNs: Fractional Physics-Informed Neural Networks},
journal = {SIAM Journal on Scientific Computing},
volume = {41},
number = {4},
pages = {A2603-A2626},
year = {2019},
doi = {10.1137/18M1229845},
URL = {     
        https://doi.org/10.1137/18M1229845
},
eprint = {     
        https://doi.org/10.1137/18M1229845

}

}

@article{Hu_2022,
	doi = {10.1137/21m1447039},  
	url = {https://doi.org/10.1137%2F21m1447039},  
	year = 2022,
	month = {sep},  
	publisher = {Society for Industrial {\&} Applied Mathematics ({SIAM})},  
	volume = {44},  
	number = {5},  
	pages = {A3158--A3182},  
	author = {Zheyuan Hu and Ameya D. Jagtap and George Em Karniadakis and Kenji Kawaguchi},  
	title = {When Do Extended Physics-Informed Neural Networks ({XPINNs}) Improve Generalization?},  
	journal = {{SIAM} Journal on Scientific Computing}
}

@article{L2_HJB,
  title={Is $ {L}^2$ Physics Informed Loss Always Suitable for Training Physics Informed Neural Network?},
  author={Wang, Chuwei and Li, Shanda and He, Di and Wang, Liwei},
  journal={Advances in Neural Information Processing Systems},
  volume={35},
  pages={8278--8290},
  year={2022}
}

@article{Ming_2021,
	doi = {10.4208/cicp.oa-2020-0219},
  
	url = {https://doi.org/10.4208%2Fcicp.oa-2020-0219},
	year = 2021,
	month = {jun},
	publisher = {Global Science Press},
	volume = {29},
  
	number = {5},
  
	pages = {1365--1384},  
	author = {Liao, Yulei and Ming, Pingbing},
	title = {Deep Nitsche Method: Deep Ritz Method with Essential Boundary Conditions},  
	journal = {Communications in Computational Physics}
}

@article{Minakowski_2023,
	doi = {10.1016/j.cam.2022.114845},  
	url = {https://doi.org/10.1016%2Fj.cam.2022.114845},  
	year = 2023,
	month = {mar},  
	publisher = {Elsevier {BV}
},  
	volume = {421},  
	pages = {114845},  
	author = {P. Minakowski and T. Richter},  
	title = {A priori and a posteriori error estimates for the Deep Ritz method applied to the Laplace and Stokes problem},
  
	journal = {Journal of Computational and Applied Mathematics}
}

@ARTICLE{shin_2020,
	author  = {Yeonjong Shin and Zhongqiang  Zhang and George Em Karniadakis},
	title   = {Error Estimates of Residual Minimization using Neural Networks for Linear PDEs},
	journal = {Journal of Machine Learning for Modeling and Computing},
	issn    = {2689-3967},
	year    = {2023},
	volume  = {4},
	number  = {4},
	pages   = {73--101}
}

@inproceedings{deRyck_2022,
 author = {De Ryck, Tim and Mishra, Siddhartha},
 booktitle = {Advances in Neural Information Processing Systems},
 editor_IGNORE = {S. Koyejo and S. Mohamed and A. Agarwal and D. Belgrave and K. Cho and A. Oh},
 pages = {10945--10958},
 publisher = {Curran Associates, Inc.},
 title = {Generic bounds on the approximation error for physics-informed (and) operator learning},
 url = {https://proceedings.neurips.cc/paper_files/paper/2022/file/46f0114c06524debc60ef2a72769f7a9-Paper-Conference.pdf},
 volume = {35},
 year = {2022}
}

@article{wang_gradientflow,
author = {Wang, Sifan and Teng, Yujun and Perdikaris, Paris},
title = {Understanding and Mitigating Gradient Flow Pathologies in Physics-Informed Neural Networks},
journal = {SIAM Journal on Scientific Computing},
volume = {43},
number = {5},
pages = {A3055-A3081},
year = {2021},
doi = {10.1137/20M1318043},
URL = {             https://doi.org/10.1137/20M1318043
},
eprint = {     
        https://doi.org/10.1137/20M1318043
}
,
}

@article{WANG_NTK,
title = {When and why PINNs fail to train: A neural tangent kernel perspective},
journal = {Journal of Computational Physics},
volume = {449},
pages = {110768},
year = {2022},
issn = {0021-9991},
doi = {https://doi.org/10.1016/j.jcp.2021.110768},
url = {https://www.sciencedirect.com/science/article/pii/S002199912100663X},
author = {Sifan Wang and Xinling Yu and Paris Perdikaris},
}

@Article{Fjordholm2017,
author={Fjordholm, U. S.
and Lanthaler, S.
and Mishra, S.},
title={Statistical Solutions of Hyperbolic Conservation Laws: Foundations},
journal={Archive for Rational Mechanics and Analysis},
year={2017},
month={Nov},
day={01},
volume={226},
number={2},
pages={809-849},
issn={1432-0673},
doi={10.1007/s00205-017-1145-9},
url={https://doi.org/10.1007/s00205-017-1145-9}
}

@article{PATEL2022110754,
title = {Thermodynamically consistent physics-informed neural networks for hyperbolic systems},
journal = {Journal of Computational Physics},
volume = {449},
pages = {110754},
year = {2022},
issn = {0021-9991},
doi = {https://doi.org/10.1016/j.jcp.2021.110754},
url = {https://www.sciencedirect.com/science/article/pii/S0021999121006495},
author = {Ravi G. Patel and Indu Manickam and Nathaniel A. Trask and Mitchell A. Wood and Myoungkyu Lee and Ignacio Tomas and Eric C. Cyr}
}

@book{bochev2009least,
  title={Least-Squares Finite Element Methods},
  author={Bochev, P.B. and Gunzburger, M.D.},
  isbn={9780387689227},
  lccn={2008943966},
  series={Applied Mathematical Sciences},
  url={https://books.google.de/books?id=5zE\_XIl-fNQC},
  year={2009},
  publisher={Springer New York}
}

@Article{Guermond2008,
author={Guermond, Jean-Luc
and Popov, Bojan},
title={$L^1$-minimization methods for Hamilton--Jacobi equations: the one-dimensional case},
journal={Numerische Mathematik},
year={2008},
month={Apr},
day={01},
volume={109},
number={2},
pages={269-284},
issn={0945-3245},
doi={10.1007/s00211-008-0142-1},
url={https://doi.org/10.1007/s00211-008-0142-1}
}

@article{Guermond_Lp,
author = {Guermond, Jean-Luc},
title = {A Finite Element Technique for Solving First-Order PDEs in $L^p$},
journal = {SIAM Journal on Numerical Analysis},
volume = {42},
number = {2},
pages = {714-737},
year = {2004},
doi = {10.1137/S0036142902417054},
URL = {   
        https://doi.org/10.1137/S0036142902417054
  },
eprint = { 
        https://doi.org/10.1137/S0036142902417054
}
}

@article{Lavery_burgers,
title = {Nonoscillatory solution of the steady-state inviscid burgers' equation by mathematical programming},
journal = {Journal of Computational Physics},
volume = {79},
number = {2},
pages = {436-448},
year = {1988},
issn = {0021-9991},
doi = {https://doi.org/10.1016/0021-9991(88)90024-1},
url = {https://www.sciencedirect.com/science/article/pii/0021999188900241},
author = {John E Lavery}
}

@article{SOD19781,
title = {A survey of several finite difference methods for systems of nonlinear hyperbolic conservation laws},
journal = {Journal of Computational Physics},
volume = {27},
number = {1},
pages = {1-31},
year = {1978},
issn = {0021-9991},
doi = {https://doi.org/10.1016/0021-9991(78)90023-2},
url = {https://www.sciencedirect.com/science/article/pii/0021999178900232},
author = {Gary A Sod}
}

@Article{Svaerd2016,
author={Sv{\"a}rd, Magnus},
title={Entropy solutions of the compressible Euler equations},
journal={BIT Numerical Mathematics},
year={2016},
month={Dec},
day={01},
volume={56},
number={4},
pages={1479-1496},
issn={1572-9125},
doi={10.1007/s10543-016-0611-3},
url={https://doi.org/10.1007/s10543-016-0611-3}
}

@inproceedings{gemp2019,
  title={The unreasonable effectiveness of adam on cycles},
  author={Gemp, Ian and McWilliams, Brian},
  booktitle = {33rd Conference on Neural Information Processing Systems (NeurIPS 2019), Vancouver, Canada},
  year      = {2019},
  url       = {https://sgo-workshop.github.io/CameraReady2019/11.pdf},
}

@article {Dinca_pLap,
    AUTHOR = {Dinca, G. and Jebelean, P. and Mawhin, J.},
     TITLE = {Variational and topological methods for {D}irichlet problems
              with {$p$}-{L}aplacian},
   JOURNAL = {Port. Math. (N.S.)},
  FJOURNAL = {Portugaliae Mathematica. Nova S\'{e}rie},
    VOLUME = {58},
      YEAR = {2001},
    NUMBER = {3},
     PAGES = {339--378},
      ISSN = {0032-5155},
   MRCLASS = {35J65 (35J20 35J60 47H11 58E05 58E30)},
  MRNUMBER = {1856715},
MRREVIEWER = {Andrzej Szulkin},
}

@article{STRELOW2023112041,
title = {Physics informed neural networks: A case study for gas transport problems},
journal = {Journal of Computational Physics},
volume = {481},
pages = {112041},
year = {2023},
issn = {0021-9991},
doi = {https://doi.org/10.1016/j.jcp.2023.112041},
url = {https://www.sciencedirect.com/science/article/pii/S0021999123001365},
author = {Erik Laurin Strelow and Alf Gerisch and Jens Lang and Marc E. Pfetsch},
}

@article{wang2022causality,
  title={Respecting causality is all you need for training physics-informed neural networks},
  author={Wang, Sifan and Sankaran, Shyam and Perdikaris, Paris},
  journal={arXiv preprint arXiv:2203.07404},
  year={2022}
}
\appendix
\section{Time-dependent approximations with fixed zero}
\label{apx:time_dependent}
We extend the computations from section \ref{sec:hypCons} to another, somewhat more technical but also more general example.
In this example we consider another class of time-dependent approximate solutions $\tilde\soln$ with
\begin{equation}
\begin{cases}
\begin{aligned}
&\tilde\soln(t,x) \in [1-\epsilon,1+\epsilon]&& \text{for } x<-\epsilon, \\
&\tilde\soln(t,x) \in [-1-\epsilon,-1+\epsilon]&& \text{for } x\geq \epsilon,\\
&\tilde\soln(t,x) \in [1 + \epsilon, -1 - \epsilon ] && \text{for } -\epsilon \leq x \leq \epsilon
\end{aligned}
\end{cases}
\label{eqn:shockApprox2}
\end{equation}
such that there exists an  $\bar x$, independent of $t$, with $\tilde\soln(t,\bar x) = 0$. The existence of an $\bar x(t)$ follows from continuity, however the restriction to time-independent $\bar x$ is technical in nature. Further, we assume that $(\bar x - x)\tilde u(x) \geq 0 $ for all $x \in \domain$. Splitting the integration over $x$ at $\bar x$, we compute the $L^1$-Norm of the residual:
\begin{equation}
\label{eqn:l2time1}
\begin{aligned}
\norm{\residual[\tilde\soln]}_{L^1(\tDom\times(-\epsilon,\epsilon))} &= \tInt \int_{-\epsilon}^{\epsilon} \abs{f(\tilde\soln)_x + \tilde\soln_t}\dd{x}\dd{t}\\
&\geq \tInt \abs{\int_{-\epsilon}^{\bar x} f(\tilde\soln)_x+ \tilde\soln_t\dd{x}}\dd{t} +  \abs{\int_{\bar x}^{\epsilon} f(\tilde\soln)_x+ \tilde\soln_t\dd{x}}\dd{t}\\
&= \tInt \abs{f\left(\tilde\soln(t,-\epsilon)\right)- \int_{-\epsilon}^{\bar x} \tilde \soln_t \dd{x}} + \abs{f\left(\tilde\soln(t,\epsilon)\right) + \int_{\bar x}^{\epsilon} \tilde \soln_t \dd{x}}\dd{t}\\
&\geq%
\!\begin{aligned}[t]
\tInt \biggl(&f\left(\tilde\soln(t,-\epsilon)\right)+ \int_{-\epsilon}^{\bar x} \left(\tilde \soln_t\right)^{\ominus} \dd{x} - \int_{-\epsilon}^{\bar x} \left(\tilde \soln_t\right)^\oplus \dd{x}\\ &+ f\left(\tilde\soln(t,\epsilon)\right) - \int_{\bar x}^{\epsilon} \left(\tilde \soln_t\right)^\ominus \dd{x} + \int_{\bar x}^{\epsilon} \left(\tilde \soln_t\right)^\oplus \dd{x}\biggr)\dd{t}.
\end{aligned}
\end{aligned}%
\end{equation}
where $(\cdot)^\ominus$ and $(\cdot)^\oplus$ denote the negative and positive part of a function respectively.
Then, because $\bar x$ is independent of $t$, we may swap the order of integration and use the fact that
\begin{equation}
\abs{\int_0^T \left(\tilde \soln_t\right)^\oplus(t,x) - \left(\tilde \soln_t\right)^\ominus(t,x)\dd{t}} =  \abs{\int_0^T \tilde\soln_t(t,x)\dd{t}} = \abs{\tilde\soln(T,x) - \tilde\soln(0,x)} \leq 1 + \epsilon.
\end{equation}
Inserting this into the last line of \eqref{eqn:l2time1} gives 
\begin{equation}
\begin{aligned}
\norm{\residual[\tilde\soln]}_{L^1(\tDom\times(-\epsilon,\epsilon))} &\geq T(1 - \epsilon)^2 - \int_{-\epsilon}^{\bar x} 1 + \epsilon \dd{x} - \int_{\bar x}^{\epsilon} 1 + \epsilon \dd{x}\\
&\geq  T(1 - \epsilon)^2 - 2\epsilon (1 + \epsilon),
\end{aligned}
\end{equation}
and we may conclude our estimate analogously to equation \eqref{eqn:badScaling}, giving a lower bound for the $L^2$-in-space-time norm of the residual scaling as $\tfrac 1{\sqrt{\epsilon}}$.
\section{Further Implementation Notes}
\label{apx:implementation}
To improve transparency and reproducibility of our results, we outline some common and useful heuristics training generative adversarial networks that we employed in our numerical experiments. 
\subsection{Gradient Descent Algorithms}
In the general Algorithm \ref{alg:GAN}, we describe regular gradient descent weight updates. Instead, any other gradient-based optimization algorithm may be used. We have experienced good results using the Adam algorithm \cite{kingma_Adam} and its AMSGrad variant \cite{reddi_amsgrad}, which are both extremely popular choices for training neural networks due to their oftentimes faster minimization of the training loss.

\subsection{Loss Computation}
While we think of $\Loss(\soln_\theta,\testSoln_\chi,\testEnt_\nu)$ as the loss function used for minimization with respect to $\soln_\theta$ and maximization with respect to $\testSoln_\chi$ and $\testEnt_\nu$, one does not need to compute the entire loss for each of the optimization steps. As discussed during the definition of our loss \eqref{eqn:Loss}, it possesses an additive decomposition into several parts depending only on one or two networks. We do not evaluate parts of the loss when they do not contribute to the gradients for weight updates, avoiding unnecessary backpropagations.  
\subsection{Learning Rates- and Scheduling}

The learning rates $\tau_{\mathrm{min}}$ and $\tau_{\mathrm{max}}$ play an important role in training and choosing them adequately has a big impact on both training time and resulting network parameters.
On the one hand, larger learning rates are desirable because larger gradient descent steps require less overall steps, but on the other hand, because the loss function is typically non-convex, using a large and fixed learning rate often results in worse final network configurations even when the update is stable.

Instead of fixed learning rates, it is common to employ  \textit{learning rate scheduling}. It is typical to start with a large initial learning rate and then reduce it over the course of training, commonly by one or more orders of magnitude as training progresses. Further, one may \textit{cycle} learning rates in a schedule \cite{smith_cyclical}, that is, intermittently increase the learning rate again to attempt to escape local minima or speed up training when the network is in a “plateau” of the loss landscape, which is slow to traverse with small learning rates.

In our experience, cyclical learning rate scheduling was very effective at obtaining more precise results. However, our numerical experiments only use a basic linear learning rate schedule to match the original wPINNs \cite{wPINNS}, when comparing performance.

Note that it is generally unclear whether both the estimation networks $\testSoln_{\chi}$ and $\testEnt_\nu$ should use the same learning rate $\tau_{\mathrm{max}}$. If learning proves difficult for any of these networks one should consider using separate learning rates and number of maximization steps for either.

\subsection{Network Checkpointing}
An important technique during the training of neural networks is \emph{checkpointing}, i.e. saving well-performing neural network weights and returning these instead of the final weights obtained after all training epochs are completed.

As a heuristic for determining network performance, one option is the neural network loss. 
We do not choose the loss directly but rather track an exponential average $\dLoss_\mathrm{avg}$ decaying at a rate of $1-\gamma$, where $\gamma \in (0,1)$ is some parameter one is free to choose.
Updating the performance indicator $\dLoss_{\mathrm{avg}}$ is done using the loss from before the weight update to $\soln_\theta$, however typical update steps only change this quantity slightly per epoch and doing so avoids one additional loss computation per epoch, making this convenient in practice.

The exponential averaging is non-standard, however we believe it is more sensible in adversarial settings than just tracking the epoch loss directly. 
Indeed, the solution network $\soln_\theta$ may update its weights in a fashion that confuses the adversarial networks, such that they grossly underestimate the correct value of the loss temporarily.
This leads to several epochs where the  loss is low and the solution network then performs a few bad weight updates before the adversarial networks recover and give accurate estimations again. Afterwards, one correctly observes an increased loss due to the poor weight updates of previous iterations.

Thus, it is disadvantageous to save single epochs where the loss is lowest, because this favors saving networks where the adversarial networks are not working properly. Tracking an exponential average of the loss mitigates this issue, leaving only the desired property of the averaged loss decreasing as the approximation to the solution improves, when the adversarial networks are working well.
\end{document}